\documentclass[12pt]{amsart}

\usepackage{amsmath}
\usepackage{amscd}
\usepackage{amssymb}
\usepackage{latexsym}
\usepackage{stmaryrd} 
\usepackage{mathbbol}

\input xy
\xyoption{all} \CompileMatrices

\newtheorem{theorem}{Theorem}[section]
\newtheorem*{theorem*}{Theorem}
\newtheorem{lemma}[theorem]{Lemma}
\newtheorem*{lemma*}{Lemma}
\newtheorem{corollary}[theorem]{Corollary}
\newtheorem{proposition}[theorem]{Proposition}

\newtheorem{remark}[theorem]{Remark}
\newtheorem{definition}[theorem]{Definition}

\makeatletter
\def\revddots{\mathinner{\mkern1mu\raise\p@
\vbox{\kern7\p@\hbox{.}}\mkern2mu
\raise4\p@\hbox{.}\mkern2mu\raise7\p@\hbox{.}\mkern1mu}}
\makeatother 
\newcommand{\bgl}{\begin{equation}} 
\newcommand{\egl}{\end{equation}}
\newcommand{\bgloz}{\begin{equation*}} 
\newcommand{\egloz}{\end{equation*}}
\newcommand{\bgln}{\begin{eqnarray}} 
\newcommand{\egln}{\end{eqnarray}}
\newcommand{\bglnoz}{\begin{eqnarray*}} 
\newcommand{\eglnoz}{\end{eqnarray*}}
\newcommand{\btheo}{\begin{theorem}}
\newcommand{\etheo}{\end{theorem}}
\newcommand{\btheooz}{\begin{theorem*}}
\newcommand{\etheooz}{\end{theorem*}}
\newcommand{\blemma}{\begin{lemma}}
\newcommand{\elemma}{\end{lemma}}
\newcommand{\blemmaoz}{\begin{lemma*}}
\newcommand{\elemmaoz}{\end{lemma*}}
\newcommand{\bproof}{\begin{proof}}
\newcommand{\eproof}{\end{proof}}
\newcommand{\bbew}{\begin{beweis}}
\newcommand{\ebew}{\end{beweis}}
\newcommand{\bremark}{\begin{remark}\em}
\newcommand{\eremark}{\end{remark}}
\newcommand{\bdefin}{\begin{definition}}
\newcommand{\edefin}{\end{definition}}
\newcommand{\bprop}{\begin{proposition}}
\newcommand{\eprop}{\end{proposition}}
\newcommand{\bcor}{\begin{corollary}}
\newcommand{\ecor}{\end{corollary}}
\newcommand{\bfa}{\begin{cases}} 
\newcommand{\efa}{\end{cases}}
\newcommand{\cG}{\mathcal G}

\newcommand{\cL}{\mathcal L}

\newcommand{\cO}{\mathcal O}

\newcommand{\cQ}{\mathcal Q}
\newcommand{\cR}{\mathcal R}

\newcommand{\so}{\text{\tiny{$\cO$}}}

\def\Az{\mathbb{A}}
\def\Cz{\mathbb{C}}
\def\Fz{\mathbb{F}}

\def\Nz{\mathbb{N}}

\def\Qz{\mathbb{Q}}
\def\Rz{\mathbb{R}}
\def\Tz{\mathbb{T}}
\def\Zz{\mathbb{Z}}

\def\1z{\mathbb{1}}
\newcommand{\fA}{\mathfrak A}
\newcommand{\fB}{\mathfrak B}

\newcommand{\fD}{\mathfrak D}

\newcommand{\an}[1]{``#1''} 
\newcommand{\ti}{\tilde}

\newcommand{\ri}{\rightarrow}
\newcommand{\lori}{\longrightarrow}

\newcommand{\rarr}{\rightarrow}

\newcommand{\ma}{\mapsto} 
\newcommand{\loma}{\longmapsto} 
\newcommand\onto{\twoheadrightarrow} 
\newcommand\into{\hookrightarrow} 
\newcommand{\LRarr}{\Leftrightarrow} 

\newcommand{\ve}{\varepsilon}

\def\SEMI{\mbox{$\times\kern-2pt\vrule height5pt width.6pt \kern3pt $}}

\newcommand{\halb}{\tfrac{1}{2}}

\newcommand{\Aut}{{\rm Aut}\,}

\newcommand{\Spec}{{\rm Spec\,}} 

\newcommand{\id}{{\rm id}}

\newcommand{\tei}{\mid} 
\newcommand{\ntei}{\nmid} 

\newcommand{\Ad}{{\rm Ad\,}}

\newcommand{\img}{{\rm im\,}}
\renewcommand{\ker}{{\rm ker}\,}

\newcommand{\reg}{^\times} 
\newcommand{\pos}{_{>0}} 
\newcommand{\ev}{\operatorname{ev}} 
\newcommand{\abs}[1]{\lvert#1\rvert} 
\newcommand{\norm}[1]{\left\|#1\right\|} 
\newcommand{\defeq}{\mathrel{:=}} 
\newcommand{\eqdef}{\mathrel{=:}} 

\newcommand{\dop}{\text{: }} 
\newcommand{\fuer}{\text{ for }} 
\newcommand{\falls}{\text{ if }} 
\newcommand{\fa}{\text{ for all }} 
\newcommand{\ilim}{\varinjlim} 
\newcommand{\plim}{\varprojlim} 
\newcommand{\e}{{\rm e}} 
\newcommand{\Ell}{{\rm L}} 
\newcommand{\extalg}{\Lambda^* \,} 
\newcommand{\rte}{\overset{e}{\rtimes}} 
\newcommand{\lge}{\left\{} 
\newcommand{\rge}{\right\}} 
\newcommand{\lru}{\left(} 
\newcommand{\rru}{\right)} 
\newcommand{\leck}{\left[} 
\newcommand{\reck}{\right]} 
\newcommand{\lsp}{\left\langle} 
\newcommand{\rsp}{\right\rangle} 
\newcommand{\rukl}[1]{\lru #1 \rru} 
\newcommand{\eckl}[1]{\leck #1 \reck} 
\newcommand{\gekl}[1]{\lge #1 \rge} 
\newcommand{\spkl}[1]{\lsp #1 \rsp} 
\newcommand{\menge}[2]{\gekl{ #1 \dop #2 }} 
\begin{document}

\title[C*-algebras associated with integral domains]{C*-algebras associated with integral domains and crossed products by actions on adele spaces}

\author{Joachim Cuntz and Xin Li}


\thanks{\scriptsize{Research supported by the Deutsche Forschungsgemeinschaft (SFB 478).}}

\thanks{\scriptsize{The second named author is supported by the Deutsche Telekom Stiftung.}}

\begin{abstract}
We compute the K-theory for C*-algebras naturally associated with rings of integers in number fields. 

The main ingredient is a duality theorem for arbitrary global fields. It allows us to identify the crossed product arising from affine transformations on the finite adeles with the analogous crossed product algebra over the infinite adele space. 
\end{abstract}

\maketitle


\section{Introduction}

Let $R$ be a countable ring. The elements of $R$ act by addition and multiplication on $\ell^2(R)$. Denote by $\fA [R]$ the C*-algebra generated by all the corresponding operators in $\cL(\ell^2(R))$. In \cite{Cun} the first named author had studied (using a different notation) this ring C*-algebra for $R = \Zz$ and had shown that 
it has an intriguing structure. In particular, it is purely infinite simple (thus a Kirchberg algebra) and can be described as a universal C*-algebra given by generators and relations. It is also Morita equivalent to a crossed product of the algebra of functions on the finite adele space for $\Qz$ by the action of the $ax+b$-group over $\Qz$. These results were generalized in \cite{CuLi} to the case where $R$ is an integral domain with finite quotients and in \cite{Li} to general rings.

It is an obvious problem to determine the K-theory of $\fA [R]$. In \cite{Cun} the case of $R = \Zz$ was discussed and it was stated that $K_*(\fA [\Zz])$ is given as an exterior $\Zz$-algebra with one generator for each prime number in $\Zz$. A proof for this was sketched. This proof however was not complete. Moreover, in \cite{Cun} a duality theorem was stated claiming that $\fA [\Zz]$ can be described also as a crossed product of $C_0(\Rz)$ by the natural action of the $ax+b$-group over $\Qz$. Again a proof was sketched which was not complete.

In the present paper we give complete proofs for these two results generalizing them at the same time to a substantially more general setting. It turns out that the two results are related. We first prove a general duality theorem that holds for any global field $K$ and states that the crossed products $C_0(\Az_f)\rtimes P_K$ and $C_0(\Az_{\infty}) \rtimes P_K$ are isomorphic. Here $\Az_f$ and $\Az_{\infty}$ denote the spaces of finite and infinite adeles, respectively, and $P_K$ denotes the $ax+b$-group over $K$. Both crossed products describe $\fA [\so]$ up to Morita equivalence, where $\so$ is the ring of integers in $K$. We note at this point that we also determine the crossed product $C_0(\Az) \rtimes P_K$ for the full adele space $\Az = \Az_f \times \Az_{\infty}$ and show that it is Morita equivalent to the group C*-algebra $C^*(P_K)$. Moreover, we point out that in the case of number fields, the Bost-Connes system and its generalizations considered in \cite{BoCo}, \cite{CMR}, \cite{HaPa} and \cite{LLN} are carried by a natural subalgebra of $\fA[\so]$. This is explained in \cite{CuLi} in more detail. 

In a second step, we use our duality theorem to determine the K-theory for $\fA [\so]$ in the case where $\so$ is the ring of integers in a number field $K$ which contains only $1,-1$ as roots of unity. The point is that the description of $\fA [\so]$ as $C_0(\Az_{\infty}) \rtimes P_K$ is much better suited for this computation since it allows for certain homotopy arguments which do not apply in the totally disconnected space $\Az_f$. We find that the K-theory depends on the number of real embeddings of $K$: Roughly speaking, we get the exterior $\Zz$-algebra over the torsion-free part of $(K \reg, \cdot)$. But if the number of real embeddings of $K$ is even and at least $2$, we will get an additional copy of this exterior algebra with coefficients in $\Zz / 2 \Zz$. These results indicate that the K-theory of $\fA [\so]$ as such does not contain information on the class number of $K$. Therefore, one is forced to investigate finer structures in $\fA [\so]$ to find out more about the class number (compare \cite{Li}, 6.5). 

For an arbitrary number field $K$, we also determine the K-theory of the subalgebra of $\fA [\so]$ which is generated by the addition operators for elements in $\so$ and the multiplication operators coming from the torsion-free part of $K \reg$. We think of this subalgebra as $\fA [\so]$ \an{without roots of unity}.  We find that its K-theory can be described as the exterior $\Zz$-algebra over the torsion-free part of $K \reg$, with coefficients in $\Zz$ or $\Zz / 2 \Zz$ depending on the real embeddings of $K$.

The paper is structured as follows:

In Section 2, we give an overview of certain aspects of algebraic number theory which we will need. We also briefly recall the notion of ring C*-algebras. 

Then we compute the K-theory of $\fA [\Zz]$. There are several reasons why we choose to treat this special case first. On the one hand, it was this case from which all our investigations started. It serves as a guide through our computations in the general setting and thereby helps to understand the general arguments. On the other hand, at various points we can considerably shorten the calculations using special features of the case $R = \Zz$. Moreover, it is possible to work out several steps explicitly in the concrete situation of $\Zz$. This leads to a self-contained exposition, in the sense that we do not need to refer to results from algebraic number theory in this special case. 

Finally, we consider the general situation. As a first step, we establish a duality theorem for arbitrary global fields (Section 4). Using this duality, we carry out the K-theoretic computations, first for the subalgebras \an{without roots of unity} of arbitrary number fields (Section 5) and then for the whole ring C*-algebras, but under the assumption that the number fields only contain the roots of unity $1,-1$ (Section 6).

We are indebted to W. L\"{u}ck for useful discussions and for bringing Lemma \ref{BC} to our attention. 

\section{Preliminaries}

\subsection{Algebraic number theory}
\label{ANT}

Let us very briefly mention aspects from algebraic number theory which will be of interest for us. First of all, the classical objects of study in algebraic number theory are the so called number fields, which are finite (separable) extensions of $\Qz$, and the corresponding rings of integers, which are the integral closures of $\Zz$ in those fields. Moreover, it turns out that the theory of the so called function fields, which are finite separable extensions of $\Fz_p (T)$, can be - at least to some extent - developed parallely. For this reason, our constructions and some of the results (our duality theorem) will apply to both situations.

However, the final step in our K-theoretic computation is only carried out for number fields. For this, it is useful to note the following: 

\blemma
\label{pro}
Let $K$ be a number field and $\mu$ be the set of roots of unity in $K$. There exists an infinitely generated free abelian group $\Gamma \subseteq K \reg$ with $K \reg = \mu \times \Gamma$. 
\elemma

\bproof
Consider the exact sequence after Corollary (3.9) in \cite{Neu}, I, \S 4.
\eproof

Now, let $K$ be a global field, which means a number field or a function field, and let $\so$ be the integral closure of $\Zz$ or $\Fz_p [T]$ in $K$. We will be concerned with the following objects:
\bgloz
  \text{The infinite adele ring } \Az_{\infty} = \prod_{v \tei \infty} K_v,
\egloz
\bgloz
  \text{the finite adele ring } \Az_f = {\prod_{v \ntei \infty}}' K_v,
\egloz
where the restricted product is taken with respect to the maximal compact subrings $\so_v \subseteq K_v $, and
\bgloz
  \text{the full adele ring } \Az = \Az_{\infty} \times \Az_f.
\egloz
These products are taken over equivalence classes of absolute values of $K$; infinite ones for $\Az_{\infty}$, finite ones for $\Az_f$. At this point, we should note that for function fields, we choose the absolute values satisfying $\abs{T}_v > 1$ to be the infinite ones (compare \cite{Weil}, III, \S 1, Theorem 2). Instead of equivalence classes of absolute values, one can equivalently take equivalence classes of embeddings of $K$ into locally compact, nondiscrete fields, these are called places in \cite{Weil}. Thus, we can always embed $K$ diagonally into $\Az_{\infty}$, $\Az_f$ or $\Az$ as each absolute value $v$ (we choose one representative for each class) gives rise to an embedding $K \into K_v$. We will not distinguish between $K$ and its image under these embeddings. This will be our convention in general as it will become clear from the context into which object we embed.

\bremark
Adeles and their multiplicative analogues, the ideles, play an important role in class field theory. The reader is refered to \cite{Neu} or \cite{Weil} for more information.
\eremark

Starting with $\so$, we can form the profinite completion $\plim \gekl{\so / I}$ over the set of nontrivial ideals in $\so$ ordered by inclusion. It turns out that this completion coincides with the maximal compact subring $\prod_{v \ntei \infty} \so_v$ of $\Az_f$: 
\bgl
\label{profin=prod}
  \plim_{(0) \neq I} \so / I \cong \prod \so_v.
\egl
Moreover, we have $K$ as well as $\prod \so_v$ sitting inside $\Az_f$. For their intersection, we get the following

\blemma
\label{intKprofin}
$\so = K \cap (\prod_{v \ntei \infty } \so_v)$.
\elemma

\bproof
Compare \cite{Weil}, Theorem 1 of Chapter V, \S 2 for number fields. The proof for function fields is analogous using \cite{Weil}, VI.
\eproof

For an infinite place $v$ of a number field, we either have $K_v \cong \Rz$ or $K_v \cong \Cz$. In the first case $v$ is called real, and it is called a complex place otherwise. We will write $v_{\Rz}$ for a real place and $v_{\Cz}$ for a complex one. Thus, we get 
\bgloz
  \Az_{\infty} \cong \Rz^{\# \gekl{v_{\Rz}}} \times \Cz^{\# \gekl{v_{\Cz}}}
\egloz
as topological rings. Note that we consider equivalence classes of embeddings, which means that two complex embeddings which are conjugate give rise to the same place. As additive topological groups, we have $\Az_{\infty} \cong \Rz^n$ where $n$ is the degree of $K$ over $\Qz$. 

The last point we would like to talk about is duality. Let $K$ be a global field.

\btheo
\label{dual}
There exists a nontrivial character $\chi$ of $\Az$ which is trivial on $K$. Any such character yields a pairing
\bgloz
  \Az \times \Az \lori \Tz \text{; } (x,y) \loma \spkl{x,y} \defeq \chi(xy).
\egloz
This pairing induces an isomorphism of topological groups: $\Az \cong \widehat{\Az}; x \loma \eckl{y \ma \spkl{y,x}}$. Thus, we also obtain
\bgl
\label{dualK}
  \Az / K \cong \widehat{K} \text{ via } \pi(x) \loma \eckl{a \ma \spkl{a,x}},
\egl
where $\pi $ is the projection $\Az \lori \Az / K $. Moreover, we can choose $\chi$ so that $\spkl{\cdot,\cdot}$ restricted to $\Az_{\infty}$ yields an isomorphism
\bgl
\label{dualinf}
  \Az_{\infty} \cong \widehat{\Az_{\infty}} \text{; } t \loma \eckl{ s \ma \spkl{t,s}}.
\egl
\etheo

\bproof
For the first pairing, compare \cite{Weil}, IV, \S 2, Theorem 3, or \cite{Lang}, XIV, \S 6, Theorem 10. The second result is proven in \cite{Lang}, XIV, \S 1, Theorem 1 for number fields, and follows from \cite{Weil}, II, \S 5, Theorem 3 in the general case.
\eproof

\subsection{Review of our constructions}
\label{rev}

Let us recall the concept of ring C*-algebras. We will only consider the integral closures of $\Zz$ or $\Fz_p [T]$ in a global field. This is a nice situation as far as the construction of ring C*-algebras is concerned because these rings are integral domains with finite quotients. We mention that it is possible to extend the construction to arbitrary rings (see \cite{Li}). 

Now, let $\so$ be the integral closure of $\Zz$ or $\Fz_p [T]$ in a global field $K$. Consider the following operators on the Hilbert space $\ell^2 (\so)$:
\bglnoz
  && U^a \xi_r = \xi_{a+r} \fuer a \in \so, \\
  && S_b \xi_r = \xi_{br} \fuer b \in \so \reg.
\eglnoz
Here $\so \reg$ is the set of nontrivial elements in $\so$. 

In analogy to the group case, we define the reduced ring C*-algebra as 
\bgloz
  \fA_r [\so] \defeq C^*\rukl{ \menge{U^a, S_b}{a \in \so, b \in \so \reg} } \subseteq \cL(\ell^2(\so)).
\egloz

The full ring C*-algebra $\fA [\so]$ is defined as the universal C*-algebra generated by unitaries $\menge{u^a}{a \in \so}$ and isometries $\menge{s_b}{b \in \so \reg}$ satisfying
\bglnoz
  &I.& u^a s_b u^c s_d = u^{a+bc} s_{bd} \\
  &II.& \sum u^a s_b s_b^* u^{-a} = 1
\eglnoz
where we sum over $\so / (b) = \gekl{a+(b)}$ in II. 

We use the notation $(b) \defeq b \cdot \so$ for principal ideals. 

There is a canonical homomorphism $\pi: \fA [\so] \lori \fA_r [\so]$ which is called the regular representation (as for groups). Moreover, it turns out that $\fA [\so]$ is purely infinite and simple (see \cite{CuLi}, Theorem 1), so that $\pi$ is an isomorphism. This allows us to identify $\fA [\so]$ with its image under $\pi$ on $\ell^2(\so)$. 

These ring C*-algebras are closely related to the number-theoretic objects we introduced before. Namely, it turns out that there is a canonical commutative C*-subalgebra $\fD [\so] \defeq C^*( \menge{u^a e_b u^{-a}}{a \in \so, b \in \so \reg} )$, where $e_b$ is the range projection $s_b s_b^*$ of $s_b$. And the bridge to number theory is built by the observation $\Spec \fD [\so] \cong \prod \so_v$ (see \cite{CuLi}, Observation 1; the argument essentially uses (\ref{profin=prod})). Furthermore, $\fA [\so]$ can be described as a semigroup crossed product (see \cite{CuLi}, Remark 3 and \cite{La}): 
\bgl
\label{semicropro}
  \fA [\so] \cong \fD [\so] \rte \so \rtimes \so \reg \sim_M C_0(\Az_f) \rtimes K \rtimes K \reg = C_0(\Az_f) \rtimes P_K, 
\egl
where $\rte$ denotes the crossed product by endomorphisms (following the notation in \cite{Li}). Recall that, by definition, the $ax+b$-group $P_K$ is $K \rtimes K \reg$. 

From now on, we will omit the argument $\so$ and simply write $\fA$, $\fA_r$ or $\fD$ respectively as it will become clear from the context which ring we mean. 

\section{Computations for $\Qz$}
\label{Q}

As announced, we study a special case first: the integers $\Zz$ in $\Qz$. This leads to a program which serves as a guide through more general computations. Roughly speaking, the idea is to compute the K-groups explicitly for the C*-subalgebra $C^*(u^1, s_{-1}, \menge{e_b}{b \in \Zz \reg})$ by choosing a suitable filtration and then to apply the Pimsner-Voiculescu sequence iteratively to get the K-theory of the whole ring C*-algebra. Actually, this idea is already contained in \cite{Cun}. However, to work out the second step rigorously we will need as a new ingredient the comparison between finite and infinite places. 

We obtain, as announced in \cite{Cun}, $K_*(\fA) \cong \extalg (\Qz \pos)$ as $\Zz / 2 \Zz$-graded abelian groups. Here and in the sequel, $\extalg$ denotes the $\Zz / 2 \Zz$-graded exterior $\Zz$-algebra and $K_*$ is the direct sum of $K_0$ and $K_1$ with the canonical grading. 

\subsection{K-theoretic computations I}
\label{KIQ}

\blemma
\label{KA^0Q}
$K_0 (C^*(u^1, s_{-1}, \gekl{e_b})) \cong \Qz \oplus \Zz$ via 
\bglnoz
  && \eckl{e_b}_0 \ma (\tfrac{1}{b}, 0) \fa b \in \Zz \pos, \\
  && \eckl{\halb (1 + s_{-1})}_0 - \eckl{\halb (1 + u^1 s_{-1})}_0 \ma (0, 1)
\eglnoz
whereas $K_1 (C^*(u^1, s_{-1}, \gekl{e_b}))$ is trivial.

Moreover, we have the following relation in $K_0 (C^*(u^1, s_{-1}, \gekl{e_b}))$:
\bgl
  \label{KRel}
  \eckl{1}_0 = 2 \cdot \eckl{\halb (1 + u^1 s_{-1})}_0.
\egl
\elemma

We write $\eckl{\cdot}_0$ or $\eckl{\cdot}_1$ for the classes in $K_0$ or $K_1$ respectively. Furthermore, we use curly brackets $\gekl{\cdot}$ to indicate that we consider a whole family of generators of a certain type. For instance, $\gekl{e_b}$ means $\menge{e_b}{b \in \so \reg}$. 

\bproof
By universal relation II., $e_b$ lies in $C^*(u^1, e_{bd})$. Thus, we can form $\ilim \gekl{C^*(u^1, s_{-1}, e_{b_i})}$ (over $\Zz \pos$ ordered by divisibility), and we get in the inductive limit $C^*(u^1, s_{-1}, \gekl{e_b})$. Therefore, to determine $K_* (C^*(u^1, s_{-1}, \gekl{e_b}))$, we have to compute $K_* (C^*(u^1, s_{-1}, e_b))$ for single $b$ and how the inclusion $\iota_{b,bd}: C^*(u^1, s_{-1}, e_b) \into C^*(u^1, s_{-1}, e_{bd})$ acts on K-theory. 

First of all, it is well-known that $C^*(u^1, s_{-1}) \cong C^*(\Zz \rtimes (\Zz / 2 \Zz)) \cong (\Cz * \Cz)^{\sim}$. This follows by comparing the universal properties of these C*-algebras. Mutually inverse isomorphisms $C^*(u^1, s_{-1}) \rightleftharpoons (\Cz * \Cz)^{\sim}$ are given by 
\bgloz
  s_{-1} \ma 2p-1, u^1 s_{-1} \ma 2q-1 \text{ and } \halb (1 + s_{-1}) \mapsfrom p, \halb (1 + u^1 s_{-1}) \mapsfrom q,
\egloz 
where $p$ and $q$ are the canonical generators of $\Cz * \Cz$. 

The K-theory of $(\Cz * \Cz)^{\sim}$ is known, it is given by
\bgloz
  K_j ((\Cz * \Cz)^{\sim}) \cong 
  \begin{cases}
    \Zz [1]_0 \oplus \Zz [p]_0 \oplus \Zz [q]_0 \falls j=0 \\
    \gekl{0} \fuer j=1.
  \end{cases}
\egloz
This determines the K-groups of $C^*(u^1, s_{-1})$. Let us fix the identification
\bgloz
  \Zz^3 \cong K_0 (C^*(u^1, s_{-1})); 
  e_1 \ma \eckl{1}_0, e_2 \ma \eckl{\halb (1 + u^1 s_{-1})}_0, e_3 \ma \eckl{\halb (1 + s_{-1})}_0.
\egloz

This also allows us to compute $K_*(C^*(u^1, s_{-1}, e_b))$ for any $b \in \Zz \pos$ since 
\bgl
\label{matQ}
  C^*(u^1, s_{-1}, e_b) \cong M_b(C^*(u^1, s_{-1})).
\egl
The idea is that the projections $e_b$, $u^1 e_b u^{-1}$, ..., $u^{b-1} e_b u^{-(b-1)}$ decompose $\ell^2(\Zz)$ into $b$ mutually isomorphic subspaces $\ell^2(b \Zz)$, $\ell^2(1 + b \Zz)$, ..., $\ell^2((b-1)+b \Zz)$ (see Lemma \ref{mat} for more details). Thus, 
\bgl
\label{KCbQ}
  K_j (C^*(u^1, s_{-1}, e_b)) \cong 
  \begin{cases}
    \Zz^3 \falls j=0 \\
    \gekl{0} \fuer j=1.
  \end{cases}
\egl

From these calculations, it already follows that $K_1(C^*(u^1, s_{-1}, \gekl{e_b}))$ is trivial. 

It remains to compute $K_0(\iota_{b,bd})$. However, it turns out that taking (\ref{KCbQ}) into account, we get $K_0(\iota_{b,bd})=K_0(\iota_{b',b'd})$ for any $b,b' \in \Zz \pos$ (see the proof of Lemma \ref{KB^0}). Thus, it suffices to consider $\iota_d \defeq \iota_{1,d}$. Under the identification (\ref{matQ}), we get the following: 

For $d=2$, we have 
$\iota_2 (u^1) =
\rukl{
\begin{smallmatrix}
  0&u^1 \\
  1&0
\end{smallmatrix}
}$ 
and 
$\iota_2 (s_{-1}) = 
\rukl{
\begin{smallmatrix}
  s_{-1}&0 \\
  0& u^{-1} s_{-1}
\end{smallmatrix} 
}$ 
which implies on $K_0$:
\bglnoz 
  && K_0 (\iota_2)(\eckl{\halb (1 + s_{-1})}_0) 
  = \eckl{\halb (1 + s_{-1})}_0 + \eckl{\halb (1 + u^1 s_{-1})}_0 \text{ and} \\
  && K_0 (\iota_2)(\eckl{\halb (1 + u^1 s_{-1})}_0) = \eckl{1}_0.
\eglnoz 
Therefore, we get 
$K_0 (\iota_2) = 
\rukl{
\begin{smallmatrix}
  2&1&0 \\
  0&0&1 \\
  0&0&1
\end{smallmatrix}
}$. 

For $d$ odd we have 
$\iota_d (u^1) = 
\rukl{
\begin{smallmatrix}
  0 & \dotso & 0 & u^1 \\
  1 & & & 0  \\
   & \ddots & & \vdots \\
   & & 1 & 0
\end{smallmatrix}
}$ 
and
$\iota_d (s_{-1}) = 
\rukl{ 
\begin{smallmatrix}
  s_{-1} & 0 & \dotso & 0 \\
  0 & & &u^{-1} s_{-1} \\
  \vdots & & \revddots & \\
  0 &u^{-1} s_{-1}& &
\end{smallmatrix}
}$
which implies on $K_0$:
\bglnoz 
  && K_0 (\iota_d) (\eckl{\halb (1 + s_{-1})}_0) = \eckl{\halb (1 + s_{-1})}_0 + \tfrac{d-1}{2} \eckl{1}_0 \text{ and} \\ 
  && K_0 (\iota_d) (\eckl{\halb (1 + u^1 s_{-1})}_0) = \eckl{\halb (1 + u^1 s_{-1})}_0 + \tfrac{d-1}{2} \eckl{1}_0.
\eglnoz
Thus we conclude that 
$K_0 (\iota_d) = 
\rukl{
\begin{smallmatrix}
  d&\tfrac{d-1}{2}&\tfrac{d-1}{2} \\
  0&1&0 \\
  0&0&1
\end{smallmatrix}
}$. 

Putting these facts together, we get by choosing a cofinal sequence $b_i$ in $\Zz \pos$ with $b_{i+1} = 2 d_i b_i$: 
\bgl
\label{K0CQ}
  K_0 (C^*(u^1, s_{-1}, \gekl{e_b})) 
  \cong \ilim \lge \Zz^3 ;
  \lru
  \begin{smallmatrix} 
  2 d_i & d_i & d_i - 1 \\
  0 & 0 & 1   \\
  0 & 0 & 1
  \end{smallmatrix}
  \rru
  \rge
  \cong \Qz \oplus \Zz.
\egl
The map of the $i$-th $K_0$-group $\Zz^3$ into $\Qz \oplus \Zz$ is given by
\bgloz
  \Zz^3 \lori \Qz \oplus \Zz; (x, y, z) \ma (\tfrac{1}{b_i} (x + \halb y + \halb z), y).
\egloz 
This immediately implies $\eckl{1}_0 = 2 \cdot \eckl{\halb (1 + u^1 s_{-1})}_0$. Moreover, $K_0$ is generated by $\eckl{e_b}_0$ corresponding to $(\tfrac{1}{b}, 0)$ and $\eckl{\halb (1 + s_{-1})}_0 - \eckl{\halb (1 + u^1 s_{-1})}_0$ which corresponds to $(0, 1)$ under the identification in (\ref{K0CQ}). 
\eproof

The next step is to adjoin the isometries $s_b$. We consider to this end 
\bgloz
  \fA^{(m)} \defeq C^*(u^1, s_{-1}, \gekl{e_b}, s_{p_1}, \dotsc, s_{p_m}).
\egloz 
Here, $p_1 < p_2 < \dotsb$ are the prime numbers in $\Zz \pos$. By construction, we have $\fA \cong \ilim \gekl{\fA^{(0)} \into \fA^{(1)} \into \dotso}$. Therefore, it suffices to determine $K_*(\fA^{(m)})$. Similarly to (\ref{semicropro}), $\fA^{(m)}$ can be described as a semigroup crossed product. This yields 
\bgln
\label{A^mQ}
  && \fA^{(m)} \sim_M C_0(\Gamma_m \cdot (\prod \Zz_p)) \rtimes (\Gamma_m \cdot \Zz) \rtimes (\mu \times \Gamma_m) \\
  &\cong& \ilim \gekl{\fA^{(m-1)}; \Ad(s_{p_m})} \rtimes_{\Ad(s_{p_m})} \Zz \nonumber
\egln
where $\Gamma_m = \spkl{p_1, \dotsc p_m} \subseteq \Qz \reg$. We have taken the inductive limit of 
\bgloz
  \gekl{\dotso \overset{\Ad(s_{p_m})}{\lori} \fA^{(m-1)} \overset{\Ad(s_{p_m})}{\lori} \fA^{(m-1)} \overset{\Ad(s_{p_m})}{\lori} \dotso}
\egloz 
to formally invert $\Ad(s_{p_m})$. 

Just a remark on notation: When we write a product like $\Gamma_m \cdot (\prod \Zz_p)$ (or $\Gamma_m \cdot \Zz$), it means that we embed the factors into an object carrying a multiplicative structure, for instance $\Az_f$ (or $\Qz$), and take the product there. It will be clear from the context which object we mean.

(\ref{A^mQ}) is the reason why we can apply the Pimsner-Voiculescu sequence. First, we compute: 

\blemma
\label{KA^1Q}
$K_j(\fA^{(1)}) \cong \Zz$ for $j = 0, 1$.
\elemma
\bproof
First of all, it follows from Lemma \ref{KA^0Q} that $\Ad(s_2)$ induces $\halb \id_{\Qz}$ on the summand $\Qz$ of $K_0 (\fA^{(0)})$. 

To calculate $ K_0 (\Ad(s_2)) (\eckl{\halb (1 + s_{-1})}_0 - \eckl{\halb (1 + u^1 s_{-1})}_0)$, let us consider the identification 
$\fA^{(0)} \cong M_2 (\fA^{(0)})$ analogous to (\ref{matQ}) under which 
\bglnoz
  && \halb (1 + s_{-1}) \text{ corresponds to } 
  \rukl{\begin{smallmatrix} \halb (1 + s_{-1}) & 0 \\ 0 & \halb ( 1 + u^1 s_{-1} ) \end{smallmatrix}},\\
  && \halb (1 + u^1 s_{-1}) \text{ corresponds to } 
  \halb \rukl{\begin{smallmatrix} 1 & s_{-1} \\ s_{-1} & 1 \end{smallmatrix}} \sim \rukl{\begin{smallmatrix} 1 & 0 \\ 0 & 0 \end{smallmatrix}},\\
  && \Ad(s_2) (\halb (1 + s_{-1})) \text{ corresponds to } 
    \rukl{\begin{smallmatrix} \halb (1 + s_{-1}) & 0 \\ 0 & 0 \end{smallmatrix}},\\
  && \Ad(s_2) (\halb (1 + u^1 s_{-1})) \text{ corresponds to } 
    \rukl{\begin{smallmatrix} \halb (1 + u^1 s_{-1}) & 0 \\ 0 & 0 \end{smallmatrix}}.
\eglnoz
Thus, on K-theory, this isomorphism maps both 
\bgloz
  \eckl{\halb (1 + s_{-1})}_0 - \eckl{\halb (1 + u^1 s_{-1})}_0
\egloz
and
\bgloz
  K_0 (\Ad(s_2)) (\eckl{\halb (1 + s_{-1}) }_0 - \eckl{\halb (1 + u^1 s_{-1})}_0)
\egloz
to $\eckl{\halb (1 + s_{-1})}_0 - \eckl{\halb (1 + u^1 s_{-1})}_0$, where we used (\ref{KRel}).

This shows that $ K_0 (\Ad(s_2)) $ is given by 
$\rukl{\begin{smallmatrix} \halb \id_{\Qz} & 0 \\ 0 & \id_{\Zz} \end{smallmatrix}}$ on $K_0 (\fA^{(0)}) \cong \Qz \oplus \Zz$.

Hence, the Pimsner-Voiculescu sequence applied to (\ref{A^mQ}), together with Lemma \ref{KA^0Q}, gives: 
\bgloz
  \begin{CD}
  \Qz \oplus \Zz @>
  - \id_{\Qz} \oplus 0 
  >> \Qz \oplus \Zz @>>> K_0 (\fA^{(1)}) \\
  @AAA @. @VVV \\
  K_1 (\fA^{(1)}) @<<< 0 @<<< 0 
  \end{CD}
\egloz
which implies $K_j (\fA^{(1)}) \cong \Zz$ for $j = 0, 1$.
\eproof

Actually, we can go one step further and show $(\Ad(s_3))_* = \id_{K_*(\fA^{(1)})}$, but at this point, we cannot show directly 
$(\Ad(s_{p_{m+1}}))_* = \id_{K_*(\fA^{(m)})}$ in general.

\subsection{Infinite and finite places over $\Qz$}

To solve our problem given in the last section, we compare the infinite place of $\Qz$ with the finite ones. To be more precise, our goal is to prove that the crossed products arising from the $ax+b$-group $P_{\Qz}$ acting on the finite adeles $\Az_f = \Qz \reg \cdot (\prod \Zz_p)$ and on the infinite place $\Rz$ of $\Qz$ respectively are Morita equivalent. This can be written in a slightly more complicated way as
\bgloz
  C_0(\Rz) \rtimes \Qz \rtimes \Qz \reg \sim_M C_0(\Qz \reg \cdot (\prod \Zz_p)) \rtimes (\Qz \reg \cdot \Zz) \rtimes \Qz \reg.
\egloz
The point is that we actually need this result not only for $\Qz \reg$ but - more generally - for any subgroup of $\Qz \reg$ in place of the full group $\Qz \reg$. This will be proven along the way as well. 

The central idea of the proof is that the infinite place and the finite ones are connected via duality (see Lemma \ref{dualQ}). That is why we think of our result as a duality theorem.

\subsubsection{Fourier transform for $\Rz$}

Let us consider some very basic constructions (mainly to set up the notation): 

We have an action of $\Qz$ on $C_0(\Rz)$ given by translation: 
\bgloz
  \hat{\tau}: \Qz \lori \Aut (C_0(\Rz)); \hat{\tau}_a (g) (t) = g(t-a) \fa g \in C_0(\Rz), a \in \Qz, t \in \Rz.
\egloz
Moreover, the Fourier transform on $C_c (\Rz)$ is given by
\bgloz
  F_{\Rz}: C_c (\Rz) \lori C_0(\Rz); f \loma \hat{f} = \eckl{t \ma \int_{\Rz} \e (ts) f(s) ds}, 
\egloz
where we set 
\bgloz
  \e (t) \defeq \exp (2 \pi i t)
\egloz 
and identify $\Rz$ with $\widehat{\Rz}$ by $t \loma \eckl{s \ma \e (ts)}$. $F_{\Rz}$ extends to an isomorphism $F_{\Rz}: C^*(\Rz) \lori C_0(\Rz)$. 

Now, we can consider the action $\tau: \Qz \lori \Aut(C^*(\Rz))$ given by conjugating $\hat{\tau}$ by $F_{\Rz}$. By construction, $F_{\Rz}$ is a covariant isomorphism with respect to $\tau$ and $\hat{\tau}$, and it thus extends to an isomorphism $F_{\Rz}: C^*(\Rz) \rtimes_{\tau} \Qz \lori C_0(\Rz) \rtimes_{\hat{\tau}} \Qz$. To simplify the notation, we will not distinguish between covariant homomorphisms and their extensions to crossed product algebras. $\tau$ is explicitly given by $\tau_a(f)(t) = \e (-at) f(t) \fa f \in C_c (\Rz) \subseteq C^*(\Rz)$.

Furthermore, consider the action $\hat{\beta}: \Qz \reg \lori \Aut(C_0(\Rz) \rtimes_{\hat{\tau}} \Qz)$ given by 
\bgloz
  \hat{\beta}_b (g u^a) = g(b^{-1} \sqcup) u^{ab} \fa g \in C_0(\Rz), a \in \Qz.
\egloz
Again, conjugating $\hat{\beta}$ by $F_{\Rz}$ gives an action $\beta: \Qz \reg \lori \Aut(C^*(\Rz) \rtimes_{\tau} \Qz)$ such that $F_{\Rz}$ induces an isomorphism 
\bgloz
  F_{\Rz}: C^*(\Rz) \rtimes_{\tau} \Qz \rtimes_{\beta} \Qz \reg 
  \lori C_0(\Rz) \rtimes_{\hat{\tau}} \Qz \rtimes_{\hat{\beta}} \Qz \reg.
\egloz 
$\beta$ is given by $\beta_b(f u^a) = \abs{b} f(b \sqcup) u^{ab}$.

\subsubsection{Identification of crossed products}
\label{croprosQ}

From this point of departure, we will now move towards the finite adeles, and the bridge between the infinite place and the finite ones is given by the additive group of our global field $\Qz$, in the following sense: Start with the action $\lambda: \Rz \lori \Aut(C^*(\Qz))$ given by 
\bglnoz
  \lambda_t(u^a) = \e (at) u^a \fa t \in \Rz, a \in \Qz,
\eglnoz
where $C^*(\Qz)$ denotes the group C*-algebra of $(\Qz,+)$. We will show that the crossed product C*-algebras 
$C^*(\Rz) \rtimes_{\tau} \Qz \text{ and } C^*(\Qz) \rtimes_{\lambda} \Rz$ are isomorphic. 

To this end, define a linear map 
\bgloz
  \varphi: C_c(\Qz, C_c(\Rz)) \lori C^*(\Qz) \rtimes_{\lambda} \Rz; \sum_a f_a u^a \loma \eckl{t \ma \sum_a \e (at) f_a (t) u^a}.
\egloz
\blemma
\label{dense}
$\varphi$ identifies $C_c(\Qz, C_c (\Rz))$ - viewed as a *-subalgebra of $C^*(\Rz) \rtimes_{\tau} \Qz$ - with the *-subalgebra 
$C_c(\Qz \times \Rz)$ of $C^*(\Qz) \rtimes_{\lambda} \Rz$.
\elemma

\bproof
This follows by computations as in the proof of Lemma \ref{cropros}. 
\eproof

\blemma
$\varphi$ extends to an isomorphism $\varphi: C^*(\Rz) \rtimes_{\tau} \Qz \cong C^*(\Qz) \rtimes_{\lambda} \Rz$.
\elemma
\bproof
$\varphi$ extends to an isometric isomorphism $\ell^1(\Qz, \Ell^1(\Rz)) \cong \Ell^1(\Rz, \ell^1(\Qz))$, where we view $\ell^1(\Qz, \Ell^1(\Rz))$ and $\Ell^1(\Rz, \ell^1 (\Qz))$ as *-subalgebras of $C^*(\Rz) \rtimes_{\tau} \Qz$ and $ C^*(\Qz) \rtimes_{\lambda} \Rz$ respectively. Moreover, $C^*(\Rz) \rtimes_{\tau} \Qz$ is the enveloping C*-algebra of $\ell^1 (\Qz, \Ell^1 (\Rz))$ and $C^*(\Qz) \rtimes_{\lambda} \Rz$ is the enveloping C*-algebra of $\Ell^1 (\Rz, \ell^1 (\Qz))$. Thus, we indeed get an isomorphism $\varphi: C^*(\Rz) \rtimes_{\tau} \Qz \cong C^*(\Qz) \rtimes_{\lambda} \Rz$ (compare the proof of Lemma \ref{cropros} for the details).
\eproof

Once again, the $\Qz \reg$-action on $ C^*(\Rz) \rtimes_{\tau} \Qz$, conjugated by $\varphi$, yields an action 
$\alpha: \Qz \reg \lori \Aut(C^*(\Qz) \rtimes_{\lambda} \Rz)$. $\alpha$ is given by the formula
\bgloz
  \alpha_b (\eckl{ t \ma \sum_a f_a(t) u^a}) = \eckl{t \ma \sum_a \abs{b} f_a(bt) u^{ab}}
\egloz
for all $\eckl{t \ma \sum_a f_a(t) u^a} \in C_c(\Rz, \ell^1 (\Qz))$.

By construction, $\varphi$ induces an isomorphism
\bgloz
  (C^*(\Rz) \rtimes_{\tau} \Qz) \rtimes_{\beta} \Qz \reg \overset{\varphi}{\underset{\cong}{\lori}} (C^*(\Qz) \rtimes_{\lambda} \Rz) \rtimes_{\alpha} \Qz \reg.
\egloz

\subsubsection{Fourier transform for $\Qz$}

At this point, the following well-known result brings the finite adele ring or rather its maximal compact subring into the game:

\blemma
\label{dualQ}
The dual group of $\Qz$ can be identified with 
\bgloz
  Y \defeq \Rz \times_{\Zz} (\prod \Zz_p) = (\Rz \times \widehat{ \Zz }) / _{(r,z) \sim (r+1,z+1)}.
\egloz
\elemma

\bproof
We use the well-known result that $\prod \Zz_p$ can be identified with $\widehat{(\Qz / \Zz)}$ via
\bgl
  \label{iddual1}
  \prod \Zz_p \ni z \loma (\eckl{\tfrac{m}{n}} \ma \e (z(n) \cdot \tfrac{m}{n})) \in \widehat{(\Qz / \Zz)}
\egl
where we view the maximal compact subring $\prod \Zz_p$ of $\Az_f$ as the projective limit of quotients of $\Zz$ which is realized as a subspace of $\Pi_{n>0} \Zz / n \Zz$.

Now, define $Y \overset{\gamma}{\lori} \widehat{\Qz}; [r,z] \loma \eckl{\tfrac{m}{n} \ma \e ((r-z(n)) \cdot \tfrac{m}{n})}$. $\gamma$ is well-defined and continuous. Since both spaces are compact, we just have to show bijectivity to prove that $\gamma$ is a homeomorphism.

To prove surjectivity, take any $\chi \in \widehat{\Qz}$. Restricting $\chi$ to $\Zz$ yields a character of $\Zz$ which is of the form $\e (r \sqcup)$ for some $r \in \Rz$. Therefore, $\chi \cdot \e (-r \sqcup)$ has constant value $1$ on $\Zz$, hence it induces a character of $\Qz / \Zz$. In other words, there exists $z \in \prod \Zz_p$ such that $\chi (\tfrac{m}{n}) \e (-r \cdot \tfrac{m}{n}) = \e (-z(n) \cdot \tfrac{m}{n}) \fa \tfrac{m}{n} \in \Qz$ because of (\ref{iddual1}). This means $\chi = \gamma([r,z])$.

$\gamma$ is injective as well: As one immediately checks, $\gamma$ is actually a group homomorphism (where addition on $Y$ is defined componentwise). Thus, we just have to show that $\gamma$ has trivial kernel. Given $[r,z] \in \ker(\gamma)$, we have $1 \equiv \gamma ([r,z]) \vert_{\Zz} = \e (r \sqcup) \vert_{\Zz}$ which implies $r \in \Zz$. Furthermore, this shows that $r-z$ is an element in $\prod \Zz_p$ yielding the trivial character on $\Qz / \Zz$. Hence, by (\ref{iddual1}), it must be $0$, which means $[r,z] \sim [0,0]$. 
\eproof

This result can be viewed as a special case of Theorem \ref{dual}, (\ref{dualK}).

$\gamma$ can be used to identify $ C^*(\Qz)$ and $C(Y)$ via the Fourier transform given by
\bgloz
  C_c (\Qz) \overset{F_{\Qz}}{\lori} C(Y); 
  F_{\Qz}(u^{m/n})([r,z]) = \ev_{m/n} (\gamma([r,z])) = \e ((r-z(n)) \cdot \tfrac{m}{n}).
\egloz
Conjugating $\lambda: \Rz \lori \Aut(C^*(\Qz))$ by $F_{\Qz}$ yields an action $\hat{\lambda}$ on $C(Y)$ given by $\hat{\lambda}_t (f) ([r,z]) = f \circ \hat{\lambda}_t^* ([r,z])$ with $\hat{\lambda}_t^* ([r,z]) = [r+t,z]$. This follows from 
\bglnoz
  && (F_{\Qz} \circ \lambda_t(u^{m/n}))([r,z]) = F_{\Qz}(\e (\tfrac{m}{n} \cdot t) u^{m/n}) ([r,z]) \\
  &=& \e (((r+t)-z(n)) \cdot \tfrac{m}{n}) = F_{\Qz} (u^{m/n}) ([r+t,z]).
\eglnoz
Again, we get an isomorphism $C^*(\Qz) \rtimes_{\lambda} \Rz \overset{F_{\Qz}}{\underset{\cong}{\lori}} C(Y) \rtimes_{\hat{\lambda}} \Rz$. 

As the last step, we describe the action $\hat{\alpha}$ of $\Qz \reg$ on $C(Y) \rtimes_{\hat{\lambda}} \Rz$ induced by $\alpha$ conjugated by $F_{\Qz}$. For any $\Qz \reg \ni b = \tfrac{m_b}{n_b}$ ($m_b \in \Zz, n_b \in \Zz_{>0}$), consider 
\bgloz
  Y \overset{\hat{\alpha}_b^*}{\lori} Y; 
  [r,z] \loma \eckl{(r-z(n_b)) \cdot b, (z(\sqcup \cdot n_b)-z(n_b)) \cdot b}.
\egloz
Multiplication with $b = \tfrac{m_b}{n_b}$ makes sense since $z(\sqcup \cdot n_b)-z(n_b)$ is in $\prod \Zz_p$ with $z(l n_b)-z(n_b) \in n_b \Zz$ for all $l \in \Zz \pos$
and because it is independent of the representation of $b$. Moreover, the expression defining $\hat{\alpha}_b^*$ is compatible with $\sim$ so that $\hat{\alpha}_b^*$ is well-defined. Furthermore, $\hat{\alpha}_b^*$ is continuous and thus a homoemorphism since $\hat{\alpha}_b^* \circ \hat{\alpha}_{1/b}^* = \id_Y$.

Now, we claim that $\hat{\alpha}: \Qz \reg \lori \Aut(C(Y) \rtimes_{\hat{\lambda}} \Rz)$ given by $\hat{\alpha}_b = F_{\Qz} \circ \alpha_b \circ F_{\Qz}^{-1}$ is of the form $\hat{\alpha}_b (f \cdot g) = ([r,z] \ma (f \circ \hat{\alpha}_b^*) \cdot (\abs{b} g(b \sqcup))) \fa f \in C(Y), g \in C_c(\Rz)$. This follows from
\bglnoz
  && (F_{\Qz} \circ \alpha_b (\eckl{s \ma g(s) u^a})(t)) [r,z] = \abs{b} g(bt) \cdot \e ((r-z(n_a n_b))ab) \\
  &=& \abs{b} g(bt) \e (((r-z(n_b)) \cdot b-((z(\sqcup \cdot n_b)-z(n_b)) \cdot b)(n_a)) \cdot a) \\
  &=& \abs{b} g(bt) F_{\Qz}(u^a)(\eckl{(r-z(n_b)) \cdot b, (z(\sqcup \cdot n_b)-z(n_b)) \cdot b})\\
  &=& (\abs{b} g(bt) F_{\Qz}(u^a) \circ \hat{\alpha}_b^*) [r,z].
\eglnoz
\bremark
It is useful to consider the action $\kappa: \Rz \rtimes \Qz \reg \rightarrow \Aut(C(Y))$ given by $\kappa(t,b)(f) = f \circ \hat{\alpha}_b^* \circ \hat{\lambda }_t^*$ where the semidirect product is taken with respect to the action $\Qz \reg \lori \Aut(\Rz); b \loma \eckl{t \ma t/b}$. $\kappa$ is a group homomorphism since $\hat{\lambda}_t^* \circ \hat{\alpha}_b^* = \hat{\alpha}_b^* \circ \hat{\lambda}_{t/b}^*$. Using a general result on crossed products by semidirect products (compare \cite{Wil}, Proposition 3.11), one immediately deduces
\bgloz
  (C(Y) \rtimes_{\hat{\lambda}} \Rz) \rtimes_{\hat{\alpha}} \Qz \reg \cong C(Y) \rtimes_{\kappa} (\Rz \rtimes \Qz \reg).
\egloz
\eremark

\bremark
Up to now, we could just as well consider a subgroup of $\Qz \reg$ instead of the whole group. So, to sum up, we have shown that for any subgroup $\Gamma$ of $\Qz \reg$, we have an isomorphism
\bgl
  \label{finres}
  (C_0(\Rz) \rtimes_{\hat{\tau}} \Qz) \rtimes_{\hat{\beta}} \Gamma 
  \cong C(Y) \rtimes_{\kappa} (\Rz \rtimes \Gamma).
\egl
\eremark

\subsubsection{Morita equivalent crossed product C*-algebras}

\bprop
\label{MgpdQ}
The transformation groupoids associated to the action of 
\bgloz
  \Rz \rtimes \Qz \reg \text{ on } Y \text{ via } [r,z] \cdot (t,b) = \hat{\alpha}_b^* \circ \hat{\lambda}_t^*([r,z]),
\egloz
denoted by $\cG$, and of 
\bgloz
  \Qz \rtimes \Qz \reg \text{ on } \Az_f \text{ by } z \cdot (a,b) = b^{-1}(z-a),
\egloz
denoted by $\ti{\cG}$, are equivalent in the sense of \cite{MRW}.
\eprop

\bproof
We will show that both groupoids are equivalent to certain subgroupoids which we can identify.

First, consider the closed subset $\ti{N} \defeq \prod \Zz_p \subseteq \Az_f = \ti{\cG}^0$. As $\Qz \reg \cdot (\prod \Zz_p) = \Az_f$, $\ti{N}$ meets every orbit in $\ti{\cG}^0$.  Moreover, the restricted range and source maps are open (details can be found in Lemma \ref{Mgpd}). Thus, by \cite{MRW}, EXAMPLE 2.7, $\ti{\cG}$ and $\ti{\cG}^{\ti{N}}_{\ti{N}}$ are equivalent, where 
\bgloz
  \ti{\cG}^{\ti{N}}_{\ti{N}} = \menge{(z,(a,b)) \in (\prod \Zz_p) \times (\Qz \rtimes \Qz \reg)}{b(z+a) \in \prod \Zz_p}.
\egloz
As a second step, consider the closed subset $\pi(\gekl{0} \times (\prod \Zz_p)) \eqdef N$ of $Y$ where $\pi$ is the canonical projection 
$\Rz \times (\prod \Zz_p) \overset{\pi}{\lori} (\Rz \times (\prod \Zz_p)) / \Zz = Y$. $N$ meets every orbit in $Y = \cG^0$ because $\bigcup_{t \in \Rz} \hat{\lambda}_t^*(N) = Y$. Again, the restricted range and source maps are open (compare Lemma \ref{Mgpd} for the details). Thus, $\cG$ and $\cG^N_N$ are equivalent by EXAMPLE 2.7 of \cite{MRW}. 

We have $\cG^N_N = \menge{([0,z],(t,b)) \in N \times (\Rz \rtimes \Qz \reg)}{\hat{\alpha}_b^* ([t,z]) \in N}$. Now, 
\bglnoz
  && \hat{\alpha}_b^* ([t,z]) = \eckl{(t-z(n_b)) \cdot b, (z(\sqcup \cdot n_b)-z(n_b)) \cdot b} \in N \\
  &\LRarr& \eckl{(t-z) \cdot m_b}(n_b) = (t-z(n_b)) \cdot m_b \in n_b \Zz \LRarr (t-z) \cdot b \in \prod \Zz_p.  
\eglnoz
In particular, this implies $t \in \Qz$. Thus, $\ti{\cG}^{\ti{N}}_{\ti{N}}$ and $\cG^N_N$ can be identified (as in Lemma \ref{Mgpd}) via $\ti{\cG}^{\ti{N}}_{\ti{N}} \ni (z,(a,b)) \loma ([0,z],(a,b^{-1})) \in \cG^N_N$.
\eproof

If we replace $\Qz \reg$ by an arbitrary subgroup $\Gamma$ of $\Qz \reg$, we have to consider the action of $(\Gamma \cdot \Zz) \rtimes \Gamma$ on $\Gamma \cdot (\prod \Zz_p)$ and the action of $\Rz \rtimes \Gamma$ on $Y$. With these modifications, everything works out as above.

\bcor
\label{res1}
$C_0(\Gamma \cdot (\prod \Zz_p)) \rtimes (\Gamma \cdot \Zz) \rtimes \Gamma \sim_M C(Y) \rtimes_{\kappa} (\Rz \rtimes \Gamma) $ for any subgroup $\Gamma $ of $\Qz \reg $.
\ecor
\bproof
This follows from Proposition \ref{MgpdQ} (applied to $\Gamma$ instead of $\Qz \reg$) together with \cite{MRW}, THEOREM 2.8, and the well-known fact that for a transformation groupoid, the (full) groupoid C*-algebra and the corresponding (full) crossed product are isomorphic.
\eproof

\bcor
\label{MQ}
For any subgroup $\Gamma$ of $\Qz \reg$, $C_0(\Gamma \cdot (\prod \Zz_p)) \rtimes (\Gamma \cdot \Zz) \rtimes \Gamma$ and $(C_0(\Rz) \rtimes_{\hat{\tau}} \Qz) \rtimes_{\hat{\beta}} \Gamma$ are Morita equivalent.
\ecor
\bproof
This result follows by combining the last corollary with (\ref{finres}).
\eproof

\subsection{K-theoretic computations II}

Corollary \ref{MQ} enables us to continue with our computations of Section \ref{KIQ}. The crucial point is that on $\Rz$, we can work with homotopies to compute the multiplicative action of $\Qz \reg$ on K-theory.

By (\ref{semicropro}), $\fA \sim_M C_0(\Az_f) \rtimes \Qz \rtimes \Qz \reg$. Thus, by Corollary \ref{MQ}, we have to determine the K-theory of $C_0(\Rz) \rtimes \Qz \rtimes \Qz \reg$. 

As a first step, the K-theory of $C_0(\Rz) \rtimes_{\hat{\beta}_{-1}} \mu$ can be computed with the help of the split exact sequence $C_0(\Rz) \rtimes \mu \into C(\Tz) \rtimes \mu \onto C^*(\mu)$ (recall $\mu = \gekl{\pm 1}$ in this case). We get
\bgl
  \label{K-Zz0}
  K_j (C_0(\Rz) \rtimes_{\hat{\beta}_{-1}} \mu)
  \cong
  \bfa
    \Zz \falls j=0 \\
    0 \fuer j=1.
  \efa
\egl
As a next step, we have 
$K_j (C_0(\Rz) \rtimes \Qz \rtimes_{\hat{\beta}_{-1}} \mu)
\cong 
\bfa
  \Qz \oplus \Zz \fuer j=0 \\
  0 \falls j=1
\efa
$
because of Lemma \ref{KA^0Q}, (\ref{A^mQ}) for $m=0$ and Corollary \ref{MQ}. 

Similarly, Lemma \ref{KA^1Q} implies $K_j (C_0(\Rz) \rtimes \Qz \rtimes (\mu \times \Gamma_1)) \cong \Zz$ ($j=0,1$) because of (\ref{A^mQ}) for $m=1$ and Corollary \ref{MQ}. Recall that $\Gamma_m$ is $\spkl{p_1, \dotsc, p_m}$, the subgroup of $\Qz \reg$ generated by the first $m$ primes. 

The inclusion $i: C_0(\Rz) \into C_0(\Rz) \rtimes \Qz$ is covariant with respect to $\hat{\beta}$ and thus induces homomorphisms between the corresponding crossed products.

\blemma
$i: C_0(\Rz) \rtimes (\mu \times \Gamma_1) \lori C_0(\Rz) \rtimes \Qz \rtimes (\mu \times \Gamma_1)$ induces $C \cdot \id_{\Zz}$ for some $0 \neq C \in \Zz$ on both $K_0$ and $K_1$.
\elemma
\bproof
First of all, we claim that $i: C_0(\Rz) \rtimes_{\hat{\beta}_{-1}} \mu \lori C_0(\Rz) \rtimes_{\hat{\tau}} \Qz \rtimes_{\hat{\beta}_{-1}} \mu$ induces $\Zz \overset{0 \oplus (C \cdot \id)}{\lori} \Qz \oplus \Zz$ on $K_0$ for some $0 \neq C \in \Zz$.

To show this, we consider the $\hat{\beta}_{-1}$-invariant inclusion $C_0(\Rz) \into C_0(\Rz) \rtimes_{\hat{\tau}_1} \Zz$. It yields, using the Pimsner-Voiculescu sequence and its naturality, the following commutative diagram with exact rows: 
\bgloz
\begin{CD}
  ... @>>> \Zz @> 2 \id >> \Zz @>>> K_1 (C_0(\Rz) \rtimes_{\hat{\beta}_{-1}} \Zz) @>>> 0 \\
  @. @VV \cong V @VV \cong V @VV K_1(i) V @. \\
  ... @>>> \Zz @> 2 \id >> \Zz @>>> K_1 ((C_0(\Rz) \rtimes_{\hat{\tau}_1} \Zz) \rtimes_{\hat{\beta}_{-1}} \Zz) @>>> ...
\end{CD}
\egloz
Therefore, $i: C_0(\Rz) \rtimes_{\hat{\beta}_{-1}} \Zz \lori (C_0(\Rz) \rtimes_{\hat{\tau}_1} \Zz) \rtimes_{\hat{\beta}_{-1}} \Zz$ does not induce the trival map on $K_1$.

Now, by \cite{Bla}, THEOREM 10.7.1 (the sequence described therein is natural with respect to covariant homomorphisms), we get the following commutative diagram with exact rows: 
\bgloz
\begin{CD}
  K_0 (C_0(\Rz) \rtimes_{\hat{\beta}_{-1}} \mu) @>>> K_1 (C_0(\Rz) \rtimes_{\hat{\beta}_{-1}} \Zz)  @>>> 0 \\
  @V K_0(i) VV @VV K_1(i) V @. \\
  K_0 ((C_0(\Rz) \rtimes_{\hat{\tau}_1} \Zz) \rtimes_{\hat{\beta}_{-1}} \mu) @>>> K_1 ((C_0(\Rz) \rtimes_{\hat{\tau}_1} \Zz) \rtimes_{\hat{\beta}_{-1}} \Zz) @>>> ...
\end{CD}
\egloz
In the commutative square, going right and then down does not yield the trivial map, and hence, $K_0(i)$ is not trivial.

As $K_0(C_0(\Rz) \rtimes_{\hat{\beta}_{-1}} \mu) \cong \Zz$ by (\ref{K-Zz0}) and $\hat{\beta}_b \sim_h \id$ on $C_0(\Rz) \rtimes_{\hat{ \beta}_{-1}} \mu$, the nontrival image of $K_0(i)$ is fixed by $K_0(\hat{\beta}_b)$ for all $b \in \Zz \pos$. Hence it follows that 
\bgloz
  C_0(\Rz) \rtimes_{\hat{\beta}_{-1}} \mu \overset{i}{\lori} C_0(\Rz) \rtimes_{\hat{\tau}} \Qz \rtimes_{\hat{\beta}_{-1}} \mu
  \cong \ilim_{b \in \Zz \pos} \gekl{ (C_0(\Rz) \rtimes_{\hat{\tau}_1} \Zz) \rtimes \mu; \hat{\beta}_b}
\egloz
does not yield the trivial homomorphism on $K_0$, either.

Now, 
$K_j (C_0(\Rz) \rtimes \Qz \rtimes \mu)
  \cong 
  \bfa
    \Qz \oplus \Zz \falls j=0 \\
    \gekl{0} \fuer j=1
  \efa
$ 
and $K_j (C_0(\Rz) \rtimes \Qz \rtimes (\mu \times \Gamma_1)) \cong \Zz$ for $j=0,1$ as we already know. Therefore, studying the Pimsner-Voiculescu sequence and going through the possibilites yield that $K_0 (\hat{\beta}_2)$ must be of the form $\Qz \oplus \Zz \overset{? \oplus \id_{\Zz}}{\lori} \Qz \oplus \Zz$ with $? \neq \id_{\Qz}$ on $K_0 (C_0(\Rz) \rtimes \Qz \rtimes_{\hat{\beta}^{-1}} \mu)$. But we have just seen that $K_0 (i)(1)$ is fixed by $K_0 (\hat{\beta}_2)$, where $1$ is the generator of 
$\Zz \cong K_0 (C_0(\Rz) \rtimes_{\hat{\beta}_{-1}} \mu)$. Thus, $K_0 (i)(1) = (0,C)$ for some $0 \neq C \in \Zz$ 
($C$ is nontrivial as $K_0 (i) \neq 0$). This proves our claim.

Secondly, the Pimsner-Voiculescu sequence, together with its naturality, implies that the assertion of the Lemma is true. 
\eproof

\btheo
We have $K_j (\fA^{(m)}) \cong \Zz^{2^{m-1}} \fa m \in \Zz \pos$ ($j=0,1$).
\etheo

\bproof
We prove by induction on $m$ that $K_j (C_0(\Rz) \rtimes \Qz \rtimes (\mu \times \Gamma_m)) \cong \Zz^{2^{m-1}}$ for $j=0,1$ and that $C_0(\Rz) \rtimes (\mu \times \Gamma_m) \overset{i}{\lori} C_0(\Rz) \rtimes \Qz \rtimes (\mu \times \Gamma_m)$ induces 
$
\rukl{
  \begin{smallmatrix}
  C & & * \\
   & \ddots & \\
  0 & & C
  \end{smallmatrix}
}
$
on K-theory.

The case $m=1$ has just been shown in the last lemma.

Now, assume that we have proven our assertion for $m$. We have (for $j=0,1$) the following commutative diagram
\bgloz
\begin{CD}
K_j (C_0(\Rz) \rtimes (\mu \times \Gamma_m)) @> 0 >> K_j (C_0(\Rz) \rtimes (\mu \times \Gamma_m) \\
@VV i_* V @VV i_* V \\
K_j (C_0(\Rz) \rtimes \Qz \rtimes (\mu \times \Gamma_m)) @> \id - (\hat{\beta}_{p_{m+1}})_*^{-1} >> K_j (C_0(\Rz) \rtimes \Qz \rtimes (\mu \times \Gamma_m) 
\end{CD}
\egloz
As we know by induction hypothesis that $K_j ((C_0(\Rz) \rtimes \Qz \rtimes (\mu \times \Gamma_m))$ ($j=0,1$) is torsion-free and that 
$i_* = 
\rukl{
  \begin{smallmatrix}
  C & & * \\
   & \ddots & \\
  0 & & C
  \end{smallmatrix}
}
$, 
$\id - (\hat{\beta}_{p_{m+1}})_*^{-1}$ must be trivial. 

Therefore, 
\bglnoz
  && K_j (C_0(\Rz) \rtimes \Qz \rtimes (\mu \times \Gamma_{m+1})) \\
  &\cong& K_j ((C_0(\Rz) \rtimes \Qz \rtimes (\mu \times \Gamma_m)) 
  \oplus K_{j+1} ((C_0(\Rz) \rtimes \Qz) \rtimes (\mu \times \Gamma_m)) \\
  &\cong& \Zz^{2^m} 
\eglnoz
for $j=0,1$ and the inclusion $i$ induces 
$
\rukl{
  \begin{array}{ccc|ccc}
  C & & * & & & \\
   & \ddots & & & * & \\
  0 & & C & & & \\
  \hline
  & & & C & & * \\
  & 0 & & & \ddots & \\
  & & & 0 & & C \\
  \end{array}
}
$
on K-theory under this decomposition of $K_j (C_0(\Rz) \rtimes \Qz \rtimes (\mu \times \Gamma_{m+1}))$ ($j=0,1$), as we wanted to prove.

Now, the theorem follows from (\ref{A^mQ}) and Corollary \ref{MQ}.
\eproof

We can instantly derive the following consequences:
\bcor
$\Ad(s_{p_{m+1}})$ induces the identity on $K_* (\fA^{(m)})$.
\ecor

\bcor 
$K_*(\fA) \cong \extalg (\Qz \pos)$ where $K_0$ corresponds to products of even and $K_1$ corresponds to products of odd numbers of pairwise distinct primes.
\ecor

\bremark
Using analogous arguments, we can determine the K-theory of $C^*(u^1, \menge{s_b}{b \in \Zz \pos})$. This case has already been investigated in \cite{Cun}, where $C^*(u^1, \menge{s_b}{b \in \Zz \pos})$ is denoted by $\cQ_{\Nz}$. Again, the main point is that Corollary \ref{MQ} allows us to compute the multiplicative action of $\Zz \pos$ or $\Qz \pos$ on K-theory. As the final result, we get $K_*(C^*(u^1, \menge{s_b}{b \in \Zz \pos})) \cong \extalg (\Qz \pos)$ where $K_0$ corresponds to products of odd numbers, $K_1$ corresponds to products of even numbers of pairwise distinct primes.
\eremark

\bremark
\label{P}
Looking back at our explicit calculations for $\Qz$, we see the following main steps:

1. Compute the K-theory of $\fA^{(0)} = C^*( \gekl{u^a}, s_{\zeta}, \gekl{e_b} )$. Here $\zeta$ is a root of unity which generates $\mu$. 

2. Compare the finite adele ring and the infinite one.

3. Show that it is enough to consider the multiplicative action of $K \reg$ on the infinite adeles.

4. Apply the Pimsner-Voiculescu sequence iteratively, together with a homotopy argument showing that the multiplicative action of the torsion-free part of $K \reg$ is trivial on K-theory.
\eremark

\section{A duality theorem}
\label{Me}

First of all, let us concentrate on the second step of our program. We can generalize Corollary \ref{MQ} to arbitrary global fields (number fields or function fields). Our result can be viewed as a duality theorem based on the duality results of Theorem \ref{dual}. So, we prove the following

\btheo
\label{M}
Let $K$ be a global field and $\Gamma$ be a subgroup of $K \reg$. 

The C*-algebras $C_0(\Az_{\infty}) \rtimes K \rtimes \Gamma$ and $C_0(\Gamma \cdot (\prod \so_v)) \rtimes (\Gamma \cdot \so) \rtimes \Gamma$ are Morita equivalent, where the groups act via inverse affine transformations.
\etheo

With $P_K = K \rtimes K \reg$ we get as a special case ($\Gamma = K \reg$): 

\bcor
$ C_0(\Az_{\infty}) \rtimes P_K \sim_M C_0(\Az_f) \rtimes P_K$.
\ecor

As in the case of $\Qz$, this result allows us to compute the action of $K \reg$ by homotopies. But first of all, let us prove Theorem \ref{M}. We need two lemmas. 

\subsection{Crossed products by subgroups of the dual group}

\blemma
\label{cropros}
Assume that $(G,+)$ is a locally compact abelian group and that $H$ is a subgroup of $\hat{G}$. Equip $H$ with a topology such that $H$ becomes a locally compact group and $\delta_h (f) = \eckl{g \loma h(g) f(g)}$; $\ve_g (\ti{f}) = \eckl{ h \loma h(-g) \ti{f} (h)}$ extend to strongly continuous actions of $H$ on $C^*(G)$ and of $G$ on $C^*(H)$ respectively.

Then $\varphi: C_c(G \times H) \lori C_c(H \times G) \text{; } f \loma \eckl{ (h,g) \ma h(-g) f(g,h)}$ extends to an isomorphism $C^*(G) \rtimes_{\delta} H \cong C^*(H) \rtimes_{\ve} G$.
\elemma
Before we come to the proof, just note that the discrete topology on $H$ is always a possible choice. Actually, this is the case of interest for our applications.

Moreover, Lemma \ref{cropros} generalizes our result in Section \ref{croprosQ}.

\bproof
The strategy is to show that $\varphi$ is an isomorphism of *-algebras and that $\varphi$ is isometric with respect to the norms $\norm{\cdot}_{L^1(H,L^1(G))}$ and $\norm{\cdot}_{L^1(G,L^1(H))}$. Then we just have to see that $C^*(G) \rtimes H$ and $C^*(H) \rtimes G$ are the enveloping C*-algebras of $L^1(H,L^1(G))$ and $L^1(G,L^1(H))$.

The central idea is that infinitesimally, we have the relation
\bgl
\label{commrel}
  w_h v_g = h(g) v_g w_h
\egl
in both crossed products $C^*(G) \rtimes_{\delta} H$ and $C^*(H) \rtimes_{\ve} G$. Here, $v_g$ and $w_h$ are the infinitesimal generators corresponding to $G$ and $H$ respectively. 

So, as a first step, integrating (\ref{commrel}) gives $\varphi (f_1 * f_2) = \varphi (f_1) * \varphi (f_2)$. Thus, $\varphi$ is multiplicative. Moreover, a simple computation shows that $\varphi$ is involutive as well. 

Secondly, applying Fubini, we see that $\varphi$ extends to an isometric isomorphism 
\bgloz
  L^1(H,L^1(G)) \cong L^1(G,L^1(H)).
\egloz

Finally, $C^*(G) \rtimes H$ is defined as the norm closure of $L^1(H,C^*(G))$ with respect to the norm 
\bgloz
  \norm{f} = \sup \menge{\norm{\pi (f)}}{\pi \text{ nondegenerate representation of } L^1(H,C^*(G))}.
\egloz
Now, we claim that we can equally well take the norm closure of $L^1(H,L^1(G))$ in the norm 
\bgloz
  \norm{f}' = \sup \menge{\norm{\pi (f)}}{\pi \text{ nondegenerate representation of } L^1(H,L^1(G))}.
\egloz
To see this, it suffices to prove $\norm{\cdot} = \norm{\cdot}'$ on $L^1(H,L^1(G))$, since this algebra is dense in $L^1(H,C^*(G))$. So, it remains to show that any nondegenerate representation of $L^1(H,L^1(G))$ extends to a representation of $L^1(H,C^*(G))$ (which will automatically be nondegenerate, too). Now, any nondegenerate representation of $L^1(H,L^1(G))$ is the integrated form of a covariant representation. Actually, one can adapt the proof of the analogous statement for $L^1(H,C^*(G))$ (see for instance \cite{Ped}). The only thing one has to use is that $L^1(G)$ has an approximate unit. But then, the integrated form of the corresponding covariant representation defines a (nondegenerate) representation of $L^1(H,C^*(G))$ extending the original one. This shows that $C^*(G) \rtimes H$ is the enveloping C*-algebra of $L^1(H, L^1(G))$. Analogously, $C^*(H) \rtimes G$ is the enveloping C*-algebra of $L^1(G, L^1(H))$. But we already know that $\varphi$ extends to an isometric isomorphism $L^1(H, L^1(G)) \cong L^1(G, L^1(H))$. Thus $\varphi$ also extends to an isomorphism $C^*(G) \rtimes H \cong C^*(H) \rtimes G$.
\eproof

\subsection{Comparison of groupoids}

As a second step, consider the following transformation groupoids which are closely related to the C*-algebras appearing in Theorem \ref{M}: Fix a subgroup $\Gamma$ of $K \reg$. Each $b \in \Gamma$ acts on $\Az_{\infty}$ via multiplication by $b^{-1}$. This gives rise to an action of $\Gamma$ on $\Az_{\infty}$ and thus to the semidirect product $\Az_{\infty} \rtimes \Gamma$. Let $\cG$ be the groupoid associated to the right action of $\Az_{\infty} \rtimes \Gamma $ on $\Az / K$ via affine transformations (given by $\pi (x) \cdot (t,b) = \pi (b((t,0) + x))$, with the canonical projection $\pi : \Az \lori \Az / K$).

Now, let $\so$ be the integral closure of $\Zz$ in $K$ if $K$ is a number field, and the integral closure of $\Fz_p [T]$ in case $K$ is a function field of characteristic $p$. Moreover, $\prod \so_v$ is the maximal compact subring of $\Az_f$, as above. $\Gamma$ acts on $\Gamma \cdot \so$ by multiplication as well (this time, we do not take inverses) and we can form $(\Gamma \cdot \so) \rtimes \Gamma$. Denote by $\ti{\cG}$ the groupoid associated to the right action of $\Gamma \cdot \so \times \Gamma$ on $\Gamma \cdot (\prod \so_v) \subseteq \Az_f$ via inverse affine transformations ($z \cdot (a,b) = b^{-1} (z-a)$).

\blemma
\label{Mgpd}
$\cG$ and $\ti{\cG}$ are equivalent as groupoids, in the sense of \cite{MRW}.
\elemma

This is the analogue of Proposition \ref{MgpdQ}, but now in the general context.

\bproof
We will use Example 2.7 of \cite{MRW} to reduce our assertion to certain subgroupoids. The remaining groupoids will even be isomorphic.

First of all, it is shown in \cite{MRW}, Example 2.7, that a locally compact (Hausdorff) groupoid $G$ is equivalent to $G_N^N$ if $N$ is a closed subset of $G^0$ such that 
\begin{itemize}
  \item[(i)] $N$ meets every orbit in $G^0$
  \item[(ii)] the restricted range and source maps $G_N \lori G$ are open.
\end{itemize}
We would like to apply this result to $\cG$ and $\ti{\cG}$: Consider the first groupoid with $N = \pi (\gekl{0} \times (\prod \so_v)) \subseteq \Az / K = \cG^0$. $N$ is closed because $\gekl{0} \times (\prod \so_v)$ is compact and $\pi$ is continuous. 

$N$ satisfies (i) since given $x = (x_{\infty}, x_f) \in \Az$, we can find $a \in K$, $z \in \prod \so_v$ such that $a+z = x_f$ 
(we have $\Az_f = K + (\prod \so_v)$, see \cite{Weil}, IV, Lemma 2.1 and \cite{Weil}, I, Corollary 4.2). Thus, $x = (x_{\infty} - a, z) + a$ which implies 
\bgloz
  \pi (x) = \pi (0,z) \cdot ( x_{\infty} - a,1) = r(\pi (0,z), ( x_{\infty} - a,1))
\egloz
where $r$ is the range map of $\cG$. This shows that $N$ meets every orbit in $\cG^0$.

To prove (ii), note that $\cG_N = \menge{(\pi (x), (t,b)) \in \cG}{\pi (b((t,0)+x)) \in N}$ by definition. $s \vert_{\cG_N}: \cG_N \lori N \text{; } (\pi (x), (t,b)) \loma \pi (b((t,0)+x))$ is open because for any open subset $U \subseteq \cG$, $s(U \cap \cG_N) = s(U \cap s^{-1} (N)) = s(U) \cap N$ is an open subset of $N$ since the source map $s$ is open. It remains to prove that $r \vert_{\cG_N}: \cG_N \loma N \text{; } (\pi (x), (t,b)) \loma \pi (x)$ is open. 
To see this, take open sets $U \subseteq \Az / K$, $V \subseteq \Az_{\infty}$ and $b \in \Gamma$. Consider the open subset $U \times (V \times \gekl{b})$ of $\cG$. It suffices to look at open sets of this form since they form a basis for the topology of $\cG$.

Now, we have $r(U \times (V \times \gekl{b}) \cap \cG_N) = U \cap \pi ((-V) \times b^{-1} \cdot (\prod \so_v))$ because of the following reason: For any $x \in \Az$, $\pi (x) \in r(U \times (V \times \gekl{b}) \cap \cG_N)$ means that $\pi (x)$ lies in $U$ and that $(\pi (x), (t,b)) \in \cG_N$ for some $t \in V$. The second statement is equivalent to: \an{There exists $t \in V$ with $\pi (x) \in \pi ((-t,0)) + b^{-1} \cdot N$} which is again equivalent to \an{$\pi (x) \in \pi ((-V) \times b^{-1} \cdot (\prod \so_v))$}. This proves our claim.

But since $\pi$ is open and $b^{-1} \cdot (\prod \so_v)$ is open in $\Az_f$, $U \cap \pi ((-V) \times b^{-1} \cdot (\prod \so_v))$ is open in $\cG^0 = \Az / K$. Therefore, (i) and (ii) hold true and thus, $\cG$ is equivalent to $\cG_N^N = \menge{(\pi (x), (t,b)) \in \Az / K \times (\Az_{\infty} \rtimes \Gamma)}{\pi (x) \in N \text{ and } \pi (b((t,0)+x)) \in N}$.

We study $\ti{\cG}$ in a similar way: Consider the closed subset $\ti{N} \defeq \prod \so_v$ of $\Gamma \cdot (\prod \so_v) = \ti{\cG}^0$. $\ti{N}$ meets every orbit in
$\ti{G}^0$ by construction.

Moreover, $\ti{\cG}_{\ti{N}}$ is given as $\menge{(z,(a,b)) \in \ti{\cG}}{b^{-1} (z-a) \in \ti{N}}$. Let $\ti{r}$, $\ti{s}$ be the range and source maps of $\ti{\cG}$. 

$\ti{s} \vert_{\ti{\cG}_{\ti{N}}}$ is open as $\ti{s} (U \cap \ti{\cG}_{\ti{N}}) = \ti{s} (U \cap \ti{s}^{-1} (\ti{N})) = \ti{s} (U) \cap \ti{N}$ is open in $\ti{N} = \prod \so_v$ for any open subset $U \subseteq \ti{\cG}$ because $\ti{s}$ is open. And $\ti{r} \vert_{\ti{\cG}_{\ti{N}}}$ is open since given any open subset $U \subseteq \Gamma \cdot (\prod \so_v)$, $\ti{r} (U \times \gekl{(a,b)} \cap \ti{\cG}_{\ti{N}}) = U \cap (a + b (\prod \so_v))$ is again open in $\Gamma \cdot (\prod \so_v) = \ti{\cG}^0$. As above, it is sufficient to consider open subsets of this type as they form a basis for the topology of $\ti{\cG}$.

Thus, we have seen that (i) and (ii) hold. This implies that $\ti{\cG}$ is equivalent to $\ti{\cG}_{\ti{N}}^{\ti{N}} = \menge{(z(a,b)) \in \Gamma \cdot (\prod \so_v) \times (\Gamma \cdot \so \rtimes \Gamma)}{z \in \prod \so_v \text{ and } b^{-1} (z-a) \in \prod \so_v}$.

Finally, we want to show that $\Phi: \ti{\cG}_{\ti{N}}^{\ti{N}} \lori \cG_N^N; (z,(a,b)) \loma (\pi (0,z), (a,b^{-1}))$ defines an isomorphism of topological groupoids. 

First of all, $\Phi$ is well-defined as $(z,(a,b)) \in \ti{\cG}_{\ti{N}}^{\ti{N}}$ means $b^{-1} (z-a) \in \ti{N} = \prod \so_v$ and therefore, $\pi (b^{-1}((a,0)+(0,z))) = \pi (b^{-1}a,b^{-1}z) = \pi (0,b^{-1}z-b^{-1}a) \in N$. Furthermore, $\Phi$ is injective since $(\pi (0,z), (a,b^{-1})) = (\pi (0,z'), (c,d^{-1}))$ implies $a=c$, $b=d$ and $(0,z-z') \in K \LRarr z=z'$. $\Phi$ is surjective: Given $(\pi (x), (t,b))$ in $\cG_N^N$, $\pi (x) \in N$ means that there exists $z \in \prod \so_v$ with $\pi (x) = \pi (0,z)$. Moreover, we know $\pi (bt,bz) = \pi (b((t,0)+(0,z)) \in N$. This implies that there exists $z' \in \prod \so_v$ with $(bt,bz)-(0,z') = (bt,bz-z') \in K$. Thus, $bt \in K$ and $bt = bz-z'$ which yield $t \in K \cap (\Gamma \cdot (\prod \so_v)) = \Gamma \cdot (K \cap (\prod \so_v)) = \Gamma \cdot \so$ by Lemma \ref{intKprofin}. Therefore, we have found $(z,(t,b^{-1})) \in \ti{\cG}_{\ti{N}}^{\ti{N}}$ with $\Phi (z,(t,b^{-1})) = (\pi (x), (t,b))$.

Moreover, it is easy to check that $\Phi$ respects the groupoid structures as well.

As the last point, we have to check that $\Phi$ is a homeomorphism. It is immediate that $\Phi$ is continuous. To prove that $\Phi^{-1}$ is continuous, choose a sequence $(z_n,(a_n,b_n)) \in \ti{\cG}_{\ti{N}}^{\ti{N}}$ with $\lim_{n \rarr \infty} \Phi (z_n,(a_n,b_n)) = \Phi (z,(a,b)) \in \cG_N^N$ for $(z,(a,b)) \in \ti{\cG}_{\ti{N}}^{\ti{N}}$. This means $\lim_{n \rarr \infty} (\pi (0,z_n),(a_n,b_n^{-1})) = (\pi (0,z),(a,b^{-1}))$. We have to show $\lim_{n \rarr \infty} (z_n,(a_n,b_n)) = (z,(a,b)) \text{ in } \ti{\cG}_{\ti{N}}^{\ti{N}}$. As $\lim_{n \rarr \infty} b_n^{-1} = b^{-1}$ in the discrete group $\Gamma$, we conclude that $b_n = b$ for almost all $n$. Thus, we can assume without loss of generality that $b_n = b$ for all $n$. Moreover, we have for all $n$ that $a_n$ lies in $K$ and $b(z_n-a) \in \ti{N} = \prod \so_v$. With $b = \tfrac{l}{m}$ for some $l$, $m$ in $\so$, $m \neq 0$, we conclude that $a_n \in ((\prod \so_v) + b^{-1} \cdot (\prod \so_v)) \cap K \subseteq (\tfrac{1}{l} \cdot (\prod \so_v)) \cap K = \tfrac{1}{l} \cdot \so$ by Lemma \ref{intKprofin}. Since $\tfrac{1}{l} \cdot \so$ sits discretely in $\Az_{\infty}$, it follows that $a_n = a$ for almost all $n$. Finally, to see that $\lim_{n \rarr \infty} z_n = z$ in $\ti{N} = \prod \so_v$, observe that $\ti{N} = \prod \so_v \lori N$; $z \loma \pi (0,z)$ is a homeomorphism as it is a  bijective, continuous map between compact (Hausdorff) spaces. 

This shows that $\ti{\cG}_{\ti{N}}^{\ti{N}} \cong \cG_N^N$ as topological groupoids. Hence, we have shown
\bgloz
  \ti{\cG} \sim_M \ti{\cG}_{\ti{N}}^{\ti{N}} \cong \cG_N^N \sim_M \cG.
\egloz
\eproof

\subsection{End of proof}

With these two lemmas, we are ready for the 

\bproof[Proof of Theorem \ref{M}]
Start with the additive action of $K$ on $C_0(\Az_{\infty})$ given by $\hat{\tau}_a (g) = g(\sqcup -a)$ for all $g \in C_0(\Az_{\infty})$. 
Since $\Az_{\infty} \cong \widehat{\Az_{\infty}}$ by Theorem \ref{dual}, (\ref{dualinf}), Fourier transform yields 
\bgloz
  C^*(\Az_{\infty}) \cong C_0(\Az_{\infty}) \text{; } 
  C_c(\Az_{\infty})\ni f \loma \eckl{t \ma \int_{\Az_{\infty}} \spkl{s,t} f(s) ds}.
\egloz
Under this identification, $\hat{\tau}$ corresponds to the action $\tau$ on $C^*(\Az_{\infty})$ given by $\tau_a (f) = \spkl{-a, \sqcup} \cdot f$. Thus, we are precisely in the situation of Lemma \ref{cropros} with $G = \Az_{\infty}$, $H = K$ ($H$ viewed as a discrete group). Lemma \ref{cropros} yields 
\bgloz
  C_0(\Az_{\infty}) \rtimes_{\hat{\tau}} K \cong C^*(\Az_{\infty}) \rtimes_{\tau} K \cong C^*(K) \rtimes_{\lambda} \Az_{\infty} 
\egloz
with $\lambda_t (u^a) = \spkl{a,t} u^a$ for all $a \in K$, $t \in \Az_{\infty}$. Applying again Fourier transform, together with Theorem \ref{dual}, (\ref{dualK}), we end up with 
\bgloz
  C^*(K) \rtimes_{\lambda} \Az_{\infty} \cong C(\hat{K}) \rtimes \Az_{\infty} \cong C(\Az / K) \rtimes_{\hat{\lambda}} \Az_{\infty} 
\egloz
with $\hat{\lambda}_t (f) = f(\pi (t,0) + \sqcup)$. 

So, to sum up these observations, we have an isomorphism 
\bgloz
  C_0(\Az_{\infty}) \rtimes_{\hat{\tau}} K \cong C(\Az / K) \rtimes_{\hat{\lambda}} \Az_{\infty}.
\egloz
Now, let $\Gamma$ be a subgroup of $K \reg$. Under the last identification, the action of $\Gamma$ on $C_0(\Az_{\infty}) \rtimes_{\hat{\tau}} K$ given by $\hat{\beta}_b (g \cdot u^a) = g(b^{-1} \sqcup) \cdot u^{ab} \text{ for all } g \in C_0(\Az_{\infty}), a \in K, b \in K \reg$ corresponds to the following action of $\Gamma$ on $C(\Az / K) \rtimes_{\hat{\lambda}} \Az_{\infty}$: $\hat{\alpha}_b (\eckl{t \loma f_t}) = \eckl{t \loma \abs{N(b)} f_{bt} (b \sqcup)}$ where $N$ denotes the norm on $K \reg$. 

Using \cite{Wil}, Proposition 3.11, we deduce 
\bgloz
  (C_0(\Az_{\infty}) \rtimes_{\hat{\tau}} K) \rtimes_{\hat{\beta}} \Gamma 
  \cong (C(\Az / K) \rtimes_{\hat{\lambda}} \Az_{\infty}) \rtimes_{\hat{\alpha}} \Gamma 
  \cong C(\Az / K) \rtimes (\Az_{\infty} \rtimes \Gamma).
\egloz
The semidirect product $\Az_{\infty} \rtimes \Gamma$ is taken with respect to the action of $\Gamma$ on $\Az_{\infty}$ which we already had in Theorem \ref{M}, and the action of $\Az_{\infty} \rtimes \Gamma$ on $C(\Az / K)$ is given by $(t,b) \cdot f(x) = f (x \cdot (t,b)) = f(\pi (b((t,0)+x)))$. Thus, Lemma \ref{Mgpd}, combined with \cite{MRW}, Theorem 2.8, yields  
\bglnoz
  && C_0(\Az_{\infty}) \rtimes (K \rtimes \Gamma) \cong (C_0(\Az_{\infty}) \rtimes_{\hat{\tau}} K) \rtimes_{\hat{\beta}} \Gamma 
  \cong C(\Az / K) \rtimes (\Az_{\infty} \rtimes \Gamma) \\
  &\cong& C^*(\cG) \sim_M C^*(\ti{\cG}) \cong C_0(\Gamma \cdot (\prod \so_v)) \rtimes (\Gamma \cdot \so) \rtimes \Gamma. 
\eglnoz 
The first identification follows again from \cite{Wil}, Proposition 3.11. Moreover, the first and the last crossed products are given by the corresponding actions via inverse affine transformations.
\eproof

\section{Computations \an{without roots of unity}}

Let us concentrate on number fields now. Fix such a field $K$ and let $\so$ be the ring of integers in $K$. Before we turn to the whole ring C*-algebra of $\so$, let us study the C*-algebra \an{without roots of unity} first. More precisely, we fix a subgroup $\Gamma$ of $K \reg$ with $K \reg = \mu \times \Gamma$ and consider the decomposition $\so \reg = \mu \times (\Gamma \cap \so \reg)$. Note that $\mu \subseteq \so$ by definition. Then we consider $\fB = C^*( \menge{u^a, s_b}{a \in \so, b \in \Gamma \cap \so \reg})$. This C*-subalgebra is generated by all unitaries given by addition but only those isometries which come from the torsion-free part $\Gamma$ of $K \reg$. We point out that $\fB$ depends on the choice of $\Gamma$. 

The reason why we first compute the K-groups of $\fB$ is twofold: On the one hand, it is possible to carry out the calculation in complete generality, in contrast to the computation of $K_*(\fA)$ (at least up to now). On the other hand, some of the results which we prove along our way to determining $K_*(\fB)$ will enter into the calculation of $K_*(\fA)$ later on.

As far as our strategy is concerned, we more or less follow the program described in Remark \ref{P}. This means that we will first compute the K-theory of $\fB^{(0)} \defeq C^*(\gekl{u^a}, \gekl{e_b})$ via a suitable filtration and then use Theorem \ref{M} to adjoin the isometries by a homotopy argument. 

The final result is 
\bglnoz
  && K_*(\fB) \\
  &\cong& 
  \bfa
    \extalg (\Gamma) \falls \# \gekl{v_{\Rz}} = 0 \\
    \extalg (\Gamma) \falls \# \gekl{v_{\Rz}} \text{odd and } \# \menge{v_{\Rz}}{v_{\Rz}(b)<0} \text{even } \forall \: b \in \Gamma \\
    (\Zz / 2 \Zz) \otimes_{\Zz} \extalg (\Gamma) \falls \# \gekl{v_{\Rz}} \text{odd and } \exists \: b \in \Gamma: \# \menge{v_{\Rz}}{v_{\Rz}(b)<0} \text{odd} \\
    (\Zz / 2 \Zz) \otimes_{\Zz} \extalg (\Gamma) \falls \# \gekl{v_{\Rz}} \geq 2 \text{ even}
  \efa
\eglnoz
Here, $\Zz / 2 \Zz$ is trivially graded and we consider graded tensor products.

\subsection{K-theory of $\fB^{(0)}$}

As a first step, the K-theory of $\fB^{(0)}$ can be computed with the help of a suitable filtration. 

\blemma
\label{KB^0}
$ K_j (\fB^{(0)}) \cong 
  \bfa
    \Qz^{2^{n-1}} \falls j \equiv n+1 \mod 2 \\
    \Qz^{2^{n-1}-1} \oplus \Zz \falls j \equiv n \mod 2
  \efa
$.

If $\omega_1, \dotsc, \omega_n$ is a $\Zz$-basis for $\so$ and if we write $u(i) \defeq u^{\omega_i}$, then the copy of $\Zz$ will be generated by $\eckl{u(1)}_1 \times \dotsb \times \eckl{u(n)}_1 \in K_n (\fB^{(0)})$. Here, we take the product on K-theory as it is described in \cite{HiRo}, 4.7.
\elemma

\bproof
By relation II., $e_b=\sum u^{ba} e_{bd} u^{-ba}$ where we sum over $\so / (d) = \gekl{a+(d)}$. Thus, we obtain the inductive system $\gekl{C^*(\gekl{u^a}, e_b); \iota_{b,bd}}$ with the inclusions $\iota_{b,bd}: C^*(\gekl{u^a}, e_b) \into C^*(\gekl{u^a}, e_{bd})$. The associated inductive limit coincides with $C^*(\gekl{u^a}, \gekl{e_b}) = \fB^{(0)}$. Hence, we have to determine $K_*(C^*(\gekl{u^a}, e_b))$ for single $b$ as well as $(\iota_{b,bd})_*$. 

First of all, $C^*(\gekl{u^a}) \cong C^*(\so) \cong C^*(\Zz^n)$. Thus, $K_*(C^*(\gekl{u^a})) \cong \Zz^{2^n}$ and $\menge{\eckl{u(i_1)}_1 \times \dotsb \times \eckl{u(i_k)}_1}{1 \leq i_1 < \dotsb < i_k \leq n, 0 \leq k \leq n}$ is a $\Zz$-basis. For $k=0$ we get $\eckl{1}_0$. Moreover, we have 

\blemma
\label{mat}
For any $b \in \Zz \pos$, we have $C^*(\gekl{u^a}, e_b) \cong M_{b^n}(C^*(\gekl{u^a}))$.
\elemma

\bproof[Proof of Lemma \ref{mat}]
Choose a minimal system $\cR_b$ of representatives for $\so / (b)$ in $\so$, in the sense that for any $c, c'$ in $\cR_b$, $c-c' \in (b)$ implies $c=c'$. We always assume $0 \in \cR_b$. From this data - using $\ell^2(\so) \cong \ell^2(c+(b)); \xi_r \ma \xi_{c+br}$ - we construct a unitary $\ell^2(\so) \cong \oplus_{c \in \cR_b} \ell^2(c+(b)) \cong \oplus_{\cR_b} \ell^2(\so)$. Conjugation with this unitary yields an isomorphism
\bgloz
  \varphi_b: \cL(\ell^2(\so)) \cong M_{b^n}(\cL(\ell^2(\so))); T \ma (s_b^* u^{-c} T u^{c'} s_b)_{c,c'}.
\egloz

We show that $\varphi_b(C^*(\gekl{u^a}, e_b)) = M_{b^n}(C^*(\gekl{u^a}))$: The universal property of $M_{b^n}(C^*(\gekl{u^a})) \cong M_{b^n}(\Cz) \otimes C^*(\gekl{u^a})$ gives rise to a homomorphism $\phi_b: M_{b^n}(C^*(\gekl{u^a})) \ri C^*(\gekl{u^a}, e_b)$ via $\phi_b(e_{c,c'} \otimes 1) = u^c e_b u^{-c'}$; $\phi_b(1 \otimes u^a) = u^{ba}$. $\phi_b$ is surjective and we have $\varphi_b \circ \phi_b = \id_{M_{b^n}(C^*(\gekl{u^a}, e_b))}$ by construction. This implies $\varphi_b(C^*(\gekl{u^a}, e_b)) = M_{b^n}(C^*(\gekl{u^a}))$.
\eproof

Thus, $C^*(\gekl{u^a}) \overset{e_{0,0} \otimes \id}{\lori} M_{b^n}(C^*(\gekl{u^a})) \overset{\phi_b}{\cong} C^*(\gekl{u^a}, e_b)$ induces an isomorphism on K-theory. By continuity of $K_*$, we get
\bglnoz
  && K_* (\fB^{(0)}) \cong \ilim_{\Zz \pos} \gekl{K_* (C^*(\gekl{u^a}, e_b)); (\iota_{b,bd})_*} \\
  &\cong& \ilim_{\Zz \pos} \gekl{K_* (C^*(\gekl{u^a})); (e_{0,0} \otimes \id)_*^{-1} \circ (\varphi_{bd})_* \circ (\iota_{b,bd})_* \circ (\phi_b)_* \circ (e_{0,0} \otimes \id)_*}
\eglnoz
where we used $\phi_{bd}^{-1} = \varphi_{bd}$ (see Lemma \ref{mat}). It suffices to take the inductive limit over $\Zz \pos$ as $\Zz \pos$ is cofinal in $\so \reg$. 

To understand the structure maps, note that by a modified version of Lemma \ref{mat}, $\phi_{bd} \circ (e_{0,0} \otimes \id)$ can be written as the composition 
\bglnoz
  && C^*(\gekl{u^a}) \overset{e_{0,0} \otimes \id}{\lori} M_{d^n}(C^*(\gekl{u^a})) \overset{\phi_d}{\lori} C^*(\gekl{u^a}, e_d) \\
  && \overset{e_{0,0} \otimes \id}{\lori} M_{b^n}(C^*(\gekl{u^a}, e_d)) \overset{\phi_b}{\lori} C^*(\gekl{u^a}, e_{bd})
\eglnoz 
and that the diagram
\bgloz
\begin{CD}
C^*(\gekl{u^a}) @> \phi_b \circ (e_{0,0} \otimes \id) >> C^*(\gekl{u^a}, e_b) \\
@VV \iota_{1,d} V @VV \iota_{b,bd} V \\
C^*(\gekl{u^a}, e_d) @> \phi_b \circ (e_{0,0} \otimes \id) >> C^*(\gekl{u^a}, e_{bd})  
\end{CD}
\egloz
commutes (these observations follow from $\phi_b \circ (e_{0,0} \otimes \id) = \Ad(s_b)$). Thus, 
\bglnoz
  && (e_{0,0} \otimes \id)_*^{-1} \circ (\varphi_{bd})_* \circ (\iota_{b,bd})_* \circ (\phi_b)_* \circ (e_{0,0} \otimes \id)_* \\
  &=& (e_{0,0} \otimes \id)_*^{-1} \circ (\varphi_d)_* \circ (e_{0,0} \otimes \id)_*^{-1} 
  \circ (\varphi_b)_* \circ (\phi_b)_* \circ (e_{0,0} \otimes \id)_* \circ  (\iota_{1,d})_* \\
  &=& (e_{0,0} \otimes \id)_*^{-1} \circ (\varphi_d)_* \circ (\iota_{1,d})_*.
\eglnoz
Therefore, with $\iota_d \defeq \iota_{1,d}$, it remains to determine $(e_{0,0} \otimes \id)_*^{-1} \circ (\varphi_d)_* \circ (\iota_d)_*$.

Now, define $\nu_d: C^*(\gekl{u^a}) \lori C^*(\gekl{u^a}); u^a \ma u^{da} $ for $d \in \Zz \pos$. Functoriality of the K-theoretic product yields 
\bgloz
  (\nu_d)_* (\eckl{u(i_1)}_1 \times \dotsb \times \eckl{u(i_k)}_1) = d^k \eckl{u(i_1)}_1 \times \dotsb \times \eckl{u(i_k)}_1.
\egloz
The crucial observation for our purposes is that we have 
\bgl
\label{cru}
  \varphi_d \circ \iota_d \circ \nu_d = 1_{d^n} \otimes \id.
\egl
This follows from the construction. We get  
\bgloz
  (e_{0,0} \otimes \id)_*^{-1} \circ (\varphi_d)_* \circ (\iota_d)_* \circ (\nu_d)_* = (e_{0,0} \otimes \id)_*^{-1} \circ (1_{d^n} \otimes \id)_* = d^n \cdot \id_*
\egloz
and thus $(e_{0,0} \otimes \id)_*^{-1} \circ (\varphi_d)_* \circ (\iota_d)_* = d^n \cdot (\nu_d)_*^{-1}$. We conclude that
\bgln
\label{adsp}
  && (e_{0,0} \otimes \id)_*^{-1} \circ (\varphi_d)_* \circ (\iota_d)_* (\eckl{u(i_1)}_1 \times \dotsb \times \eckl{u(i_k)}_1) \\
  &=& d^{n-k} \cdot (\eckl{u(i_1)}_1 \times \dotsb \times \eckl{u(i_k)}_1). \nonumber
\egln
This allows us to calculate the K-groups of $\fB^{(0)}$, and we get the desired results.
\eproof 

\subsection{K-theory of $\fB$}
\label{KB}

We distinguish between the following cases:

1. $\# \gekl{v_{\Rz}} = 0$

By \cite{Neu}, I, Proposition (8.4), there are infinitely many primes in $\Zz \subseteq \so$ which are unramified in $\so$. Hence, since $\mu$ is finite, by taking an appropriate product of unramified primes, we find a number $p \in \Zz \pos$ in $\Gamma$ which we can extend to a $\Zz$-basis for $\Gamma$ of the form $\gekl{p, p_1, p_2, \dotsc}$.

Let $\Gamma_m$ be the subgroup $\spkl{p, p_1, \dotsc, p_m}$ of $\Gamma$ and consider $\fB^{(1)} \defeq C^*(\fB^{(0)}, s_p)$. 

Since $\Ad(s_p)_* (\eckl{u(i_1)}_1 \times \dotsb \times \eckl{u(i_k)}_1) = p^{k-n} \eckl{u(i_1)}_1 \times \dotsb \times \eckl{u(i_k)}_1$ by (\ref{adsp}) and 
\bglnoz
  && \fB^{(1)} \cong \fB^{(0)} \rte_{\Ad(s_p)} \Nz \sim_M \ilim \gekl{\fB^{(0)}; \Ad(s_p)} \rtimes_{\Ad(s_p)} \Zz \\
  &\sim_M& C_0(\Gamma_0 \cdot (\prod \so_v)) \rtimes \Gamma_0 \cdot \so \rtimes \Gamma_0
\eglnoz
(similar to (\ref{semicropro})), we conclude that $K_j (\fB^{(1)}) \cong \Zz$ for $j=0,1$. This result, Lemma \ref{KB^0} and Theorem \ref{M} show that 
\bgln
\label{Kinf0}
  && K_j (C_0(\Az_{\infty}) \rtimes K) \cong K_j (C(\prod \so_v) \rtimes \so) \cong K_j (\fB^{(0)}) \\
  &\cong&
  \bfa
    \Qz^{2^{n-1}} \fuer j \equiv n+1 \mod 2 \\
    \Qz^{2^{n-1}-1} \oplus \Zz \fuer j \equiv n \mod 2; 
  \efa
  \nonumber \\
\label{Kinf1}
  && K_j (C_0(\Az_{\infty}) \rtimes K \rtimes \Gamma_0) \cong K_j (C(\Gamma_0 \cdot (\prod \so_v)) \rtimes (\Gamma_0 \cdot \so) \rtimes \Gamma_0) \\
  &\cong& K_j (\fB^{(1)}) \cong \Zz \fuer j=0,1. \nonumber
\egln
Moreover, as the multiplicative action of $K \reg$ on $\Az_{\infty}$ is homotopic to the trivial action, the Pimsner-Voiculescu sequence yields $K_j (C_0(\Az_\infty) \rtimes \Gamma_0) \cong \Zz$ ($j=0,1$). Now, we claim $K_j (C_0(\Az_\infty) \rtimes K \rtimes \Gamma_m) \cong K_j (C_0(\Az_\infty) \rtimes \Gamma_m) \cong \Zz^{2^m}$ for $j=0,1$ and that $i_m: C_0(\Az_\infty) \rtimes \Gamma_m \lori C_0(\Az_\infty) \rtimes K \rtimes \Gamma_m$, the map induced by the $K \reg$-covariant inclusion $C_0(\Az_\infty) \into C_0(\Az_\infty) \rtimes K$, yields an injective map on K-theory. Let us prove this assertion by induction on $m$.

$m=0$: We already know $K_j (C_0(\Az_\infty) \rtimes K \rtimes \Gamma_0) \cong K_j (C_0(\Az_\infty) \rtimes \Gamma_0) \cong \Zz$ for $j=0,1$. It remains to show injectivity of $(i_0)_*$: By (\ref{Kinf0}), (\ref{Kinf1}) and the Pimsner-Voiculescu sequence, $C_0(\Az_\infty) \rtimes K \into C_0(\Az_\infty) \rtimes K \rtimes \Gamma_0$ must be nontrivial on the copy of $\Zz$ in $K_j$ ($j \equiv n \mod 2$). Moreover, $C_0(\Az_\infty) \into C_0(\Az_\infty) \rtimes \so$ is injective on K-theory by the Pimsner-Voiculescu sequence. And we have $C_0(\Az_\infty) \rtimes K \cong \ilim \gekl{C_0(\Az_\infty) \rtimes \so}$ where the structure maps fix the image of $C_0(\Az_\infty)$ in $C_0(\Az_\infty) \rtimes \so$ on K-theory. It follows that $C_0(\Az_\infty) \into C_0(\Az_\infty) \rtimes K$ is injective on K-theory and its image is fixed by $\Gamma_0$ in K-theory and thus must be contained in the copy of $\Zz$ ($j \equiv n \mod 2$). Thus, $(i_0)_*$ is injective.

Now, assume that we have proven the claim for $m$. The following diagram 
\bgloz
\begin{CD}
  K_j (C_0(\Az_\infty) \rtimes \Gamma_m) @> 0 >> K_j (C_0(\Az_\infty) \rtimes \Gamma_m) \\
  @V (i_m)_* VV @V (i_m)_* VV \\
  K_j (C_0(\Az_\infty) \rtimes K \rtimes \Gamma_m) @> \id - (\hat{\beta}_{p_{m+1}})_*^{-1} >> K_j (C_0(\Az_\infty) \rtimes K \rtimes \Gamma_m)
\end{CD}
\egloz
commutes, where we used that the multiplicative action of $K \reg$ on $C_0(\Az_\infty)$ is equivariantly homotopic to the identity. As $K_j (C_0(\Az_\infty) \rtimes K \rtimes \Gamma_m)$ ($j=0,1$) is torsion-free and $(i_m)_*$ is injective by hypothesis, it follows that $(\hat{\beta}_{p_{m+1}})_* = \id_{K_* (C_0(\Az_\infty) \rtimes K \rtimes \Gamma_m)}$. Thus, we get $K_j (C_0(\Az_\infty) \rtimes K \rtimes \Gamma_{m+1}) \cong \Zz^{2^{m+1}}$ for $j=0,1$ by the Pimsner-Voiculescu sequence. Injectivity of 
$(\iota_{m+1})_* = 
  \rukl{
  \begin{smallmatrix}
  (i_m)_* & * \\
  0 & (i_m)_*
  \end{smallmatrix}
  }
$
follows by induction hypothesis. This proves our claim. 

Therefore, 
\bgloz
  K_*(C_0(\Az_\infty) \rtimes K \rtimes \Gamma) \cong \ilim_m K_*(C_0(\Az_\infty) \rtimes K \rtimes \Gamma_m) \cong \ilim_m \extalg (\Gamma_m) \cong \extalg (\Gamma).
\egloz
Since $\fB \sim_M C_0(\Az_f) \rtimes K \rtimes \Gamma \sim_M C_0(\Az_\infty) \rtimes K \rtimes \Gamma$ by the analogue of (\ref{semicropro}) and Theorem \ref{M}, we conclude $K_*(\fB) \cong \extalg (\Gamma)$.

The remaining cases are very similar to the first one. The only difference lies in the fact that the multiplicative action does not need to be homotopic to the identity any more.

2. $\# \gekl{v_{\Rz}}$ odd and $\# \menge{v_{\Rz}}{v_{\Rz}(b)<0}$ is even for all $b \in \Gamma$

We still have that the action of $\Gamma$ is equivariantly homotopic to the identity on $C_0(\Az_{\infty})$. Thus, we can adapt the computations of the first case. Again, 
\bgloz
  K_*(\fB) \cong \extalg (\Gamma).
\egloz

3. {$\# \gekl{v_{\Rz}}$ odd and $\# \menge{v_{\Rz}}{v_{\Rz}(b)<0}$ is odd for some $b \in \Gamma$}

Choose a basis $\gekl{p, p_1, \dotsc}$ of $\Gamma$ as in the first case ($\# \gekl{v_{\Rz}} = 0$), but we additionally require $\# \menge{v_{\Rz}}{v_{\Rz}(p_1)<0}$ to be odd and $\# \menge{v_{\Rz}}{v_{\Rz}(p_i)<0}$ to be even for all $i>1$. Now, let $\Gamma_m = \spkl{p, p_1, \dotsc, p_m} \subseteq \Gamma$ as above, and let $\Gamma_m'$ be the subgroup generated by $\gekl{p, p_2, \dotsc, p_m}$ such that $\Gamma_m = \Gamma_m' \times \spkl{p_1}$. 

We have $\hat{\beta}_{p_m} \sim_h \hat{\beta}_{t_m}$ where $t_m \defeq (v(p_m)/\abs{v(p_m)})_{v \tei \infty}$. For $m \neq 1$, $\hat{\beta}_{t_m}$ is of period $2$ and induces the identity on $K_* (C_0(\Az_{\infty}))$. Thus, going through the Pimsner-Voiculescu sequences, we get that $\hat{\beta}_{t_m}$ and therefore $\hat{\beta}_{p_m}$ induces the identity on $K_* (C_0(\Az_{\infty}) \rtimes \Gamma_{m-1}')$. 

Hence, for the same reasons as above (in the first case), we get
\bgloz
  K_*(C_0(\Az_{\infty}) \rtimes K \rtimes \Gamma_m') \cong \extalg (\Gamma_m').
\egloz
But now, we have to add $p_1$. $\hat{\beta}_{p_1}$ induces $- \id$ on $K_* (C_0(\Az_{\infty}) \rtimes \Gamma_m')$ as the map is of order $2$ on K-theory but gives $- \id$ on $K_* (C_0(\Az_{\infty}))$. So, using the Pimsner-Voiculescu sequence again, we see that $\hat{\beta}_{p_1}$ induces $-\id$ on $K_*(C_0(\Az_{\infty}) \rtimes \Gamma_m')$ and thus on $K_* (C_0(\Az_{\infty}) \rtimes K \rtimes \Gamma_m')$ as well. Hence we get 
\bgloz
  K_*(C_0(\Az_{\infty}) \rtimes K \rtimes \Gamma_m) \cong (\Zz / 2\Zz) \otimes_{\Zz} \extalg (\Gamma_m')
\egloz
and finally, in the inductive limit
\bgloz
  K_*(\fB) \cong K_*(C_0(\Az_{\infty}) \rtimes K \rtimes \Gamma) \cong (\Zz / 2\Zz) \otimes_{\Zz} \extalg (\Gamma).
\egloz

4. $\# \gekl{v_{\Rz}} \geq 2$ even

Since $\Rz \cdot K = \Az_{\infty}$, there must be some $b \in \Gamma$ with $\# \menge{v_{\Rz}}{v_{\Rz}(b)<0}$ odd. Then, the same arguments as in previous case show that 
\bgloz
  K_*(\fB) \cong K_*(C_0(\Az_{\infty}) \rtimes K \rtimes \Gamma) \cong (\Zz / 2\Zz) \otimes_{\Zz} \extalg (\Gamma).
\egloz 

\section{General results for $\mu = \gekl{\pm 1}$}

As in the previous section, let $K$ be a number field of degree $n = \eckl{K:\Qz}$ and let $\so$ be the ring of integers in $K$. We follow the program of Remark \ref{P} and compute the K-groups of the ring C*-algebra associated to $\so$ under the assumption that the only roots of unity in $K$ are $\pm 1$ ($\mu = \gekl{\pm 1}$). We will explain below why we cannot treat the general case up to now. 

The final result (with $K \reg = \mu \times \Gamma$) is as follows: 
\bgloz
  K_*(\fA) \cong 
  \bfa 
    K_0(C^*(\mu)) \otimes_{\Zz} \extalg (\Gamma) \falls \# \gekl{v_{\Rz}} = 0 \\
    \extalg (\Gamma) \falls \# \gekl{v_{\Rz}} \text{ odd} \\
    \extalg (\Gamma) \oplus ((\Zz / 2 \Zz) \otimes_{\Zz} \extalg (\Gamma)) \falls \# \gekl{v_{\Rz}} \geq 2 \text{ even}.
  \efa
\egloz
Here we consider graded tensor products where $K_0(C^*(\mu))$ and $\Zz / 2 \Zz$ are trivially graded. We take the diagonal grading on the direct sum. Moreover, note that the condition \an{$\# \gekl{v_{\Rz}}$ is odd} is equivalent to \an{$n$ is odd}. 

\subsection{The K-theory of $\fA^{(0)}$}

The first step is to calculate the K-groups of $C^*(\gekl{u^a}, s_{-1})$, or - what amounts to the same - of the group C*-algebra of $\so \rtimes \mu \cong \Zz^n \rtimes \Zz / 2 \Zz$. Let $\omega_1, \dotsc, \omega_n$ be a $\Zz$-basis of $\so$ and set $u(i) = u^{\omega_i}$ for $1 \leq i \leq n$. We thank W. L\"{u}ck for pointing out the following result to us:

\blemma
\label{BC}
The K-theory of $C^*(\so \rtimes \mu)$ is given by 
\bgloz
  K_0(C^*(\so \rtimes \mu)) = G_{fin} \oplus G_{inf}
\egloz
where $G_{fin} \cong \Zz^{2^n}$ is the part coming from finite subgroups and $G_{inf} \cong \Zz^{2^{n-1}}$. $K_1(C^*(\so \rtimes \mu))$ is trivial. 

Moreover, if we identify $C^*(\so \rtimes \mu)$ with $C^*(\gekl{u^a}, s_{-1})$, we obtain the following projections whose classes in $K_0(C^*(\gekl{u^a}, s_{-1}))$ form a $\Zz$-basis for $G_{fin}$: 
\bgloz
  \halb(1+u(i_1) \dotsm u(i_k) s_{-1}) \text{ where } 1 \leq i_1 < \dotsb < i_k \leq n \text{, } 0 \leq k \leq n.
\egloz
Furthermore, the inclusion 
\bgloz
  i: C^*(\gekl{u^a}) \cong C^*(\so) \into C^*(\so \rtimes \mu) \cong C^*(\gekl{u^a}, s_{-1})
\egloz
maps $K_0(C^*(\gekl{u^a}))$ into $G_{inf}$ injectively. Thus, its cokernel is finite.
\elemma

\bproof
See \cite{ELPW}, Theorem 0.4. 
\eproof

The next step is to compute the K-theory of $\fA^{(0)} \defeq C^*(\gekl{u^a}, s_{-1}, \gekl{e_b})$. 

\blemma
\label{KA^0}
$K_0(\fA^{(0)}) \cong 
\bfa
  \Zz \oplus \Qz^{2^{n-1}} \falls n \text{ odd} \nonumber \\
  \Zz \oplus \Qz^{2^{n-1}-1} \oplus \Zz \falls n \text{ even}
\efa
$. $K_1(\fA^{(0)})$ is trivial. 
\elemma

\bproof
Again, $\fA^{(0)} \cong \ilim \gekl{C^*(\gekl{u^a}, s_{-1}, e_b)}$ and with analogous arguments as in the proof of Lemma \ref{KB^0}, we get 
\bgloz
  K_*(\fA^{(0)}) \cong \ilim_{d \in \Zz \pos} \gekl{K_*(C^*(\gekl{u^a}, s_{-1})); (e_{0,0} \otimes \id)_*^{-1} (\varphi_d)_* (\iota_d)_*}.
\egloz 
It remains to determine 
\bgloz
  \kappa_d \defeq (e_{0,0} \otimes \id)_*^{-1} (\varphi_d)_* (\iota_d)_*: K_*(C^*(\gekl{u^a}, s_{-1})) \lori K_*(C^*(\gekl{u^a}, s_{-1}))
\egloz
for $d \in \Zz \pos$. In our computations of $K_*(\fB^{(0)})$ in Lemma \ref{KB^0}, we have already seen in (\ref{adsp}) that 
\bgloz
  \kappa_d(\eckl{u(i_1)}_1 \times \dotsb \times \eckl{u(i_k)}_1) = d^{n-k} \eckl{u(i_1)}_1 \times \dotsb \times \eckl{u(i_k)}_1
\egloz
for $1 \leq i_1 < \dotsb < i_k \leq n$, $0 \leq k \leq n$. Thus, by Lemma \ref{BC}, $\kappa_d \vert_{G_{inf}}$ is given by a diagonal matrix whose entries are powers of $d$. Among these diagonal entries, $1 = d^0$ appears if and only if $n$ is even, this entry then corresponds to $(\eckl{u(1)}_1 \times \dotsb \times \eckl{u(n)}_1) \in K_n(C^*(\gekl{u^a}, s_{-1}))$.

To determine $\kappa_d \vert_{G_{fin}}$, we distinguish between two cases: 

1. $d = 2$:
We choose $\cR_2 \defeq \gekl{0, \omega_1, \omega_2, \omega_1 + \omega_2, \omega_3, \dotsc, \omega_1 + \dotsb + \omega_n}$. With this choice, and under analogous identifications as in Lemma \ref{mat}, 
\bgloz
  \varphi_2 \circ \iota_2 (s_{-1}) 
  = 
  \rukl{
  \begin{smallmatrix}
  s_{-1} & & & & & \\
   & u(1)^* s_{-1} & & & & \\
   & & u(2)^* s_{-1} & & & \\
   & & & u(1)^* u(2)^* s_{-1} & & \\
   & & & & \ddots & \\
   & & & & & u(1)^* \dotsm u(n)^* s_{-1} 
  \end{smallmatrix}
  }.
\egloz
This shows that $\kappa_2(\eckl{\halb(1+s_{-1})}_0) = \sum \eckl{\halb(1 + u(i_1) \dotsm u(i_k) s_{-1})}_0$ where the sum is taken over all $1 \leq i_1 < \dotsb < i_k \leq n$, $0 \leq k \leq n$.

Now, the remaining symmetries $u(i_1) \dotsm u(i_k) s_{-1}$ with $1 \leq i_1 < \dotsb < i_k \leq n$, $1 \leq k \leq n$ map $\ell^2(c+(2))$ bijectively into $\ell^2(c'+(2))$ for $c \neq c'$ in $\cR_2$. Thus, a typical building block in $\varphi_2 \circ \iota_2(u(i_1) \dotsm u(i_k) s_{-1})$ is of the form 
$ 
  \rukl{
  \begin{smallmatrix}
  0 & \dotso & V^* \\
  \vdots & 0 & \vdots \\
  V & \dotso & 0
  \end{smallmatrix}
  }
$
for a unitary $V$. 
As 
\bgloz
  \halb
  \rukl{
  \begin{smallmatrix}
  1 & 0 & \dotso & 0 & V^* \\
  0 & & & & 0 \\
  \vdots & & 0 & & \vdots \\
  0 & & & & 0 \\
  V & 0 & \dotso & 0 & 1
  \end{smallmatrix}
  }
  \sim 
  \rukl{
  \begin{smallmatrix}
  1 & 0 & \dotso & 0 \\
  0 & & & \\
  \vdots & & 0 & \\
  0 & & &
  \end{smallmatrix}
  },
\egloz
we get $\kappa_2(\eckl{\halb(1+u(i_1) \dotsm u(i_k))}_0) = 2^{n-1} \eckl{1}_0$. Thus, $\kappa_2$ is given by 
\bgloz
  \rukl{
  \begin{array}{cccc|ccc}
  1 & & & & & & \\
  \vdots & & 0 & & & 0 & \\
  1 & & & & & & \\
  \hline
  0 & 2^{n-1} & \dotso & 2^{n-1} & 2^n & & \\
  \vdots & & 0 & & & \ddots & \\
  0 & & & & & & 2^{?}
  \end{array}
}
\egloz
where $?$ is $0$ or $1$ depending on the parity of $n$.

2. Let $d$ be odd, say $d = 2d'+1$:

This time let $\cR_d$ be $\menge{\sum_{m=1}^n l_m \omega_m}{-d' \leq l_m \leq d'}$. 

For any $0 \leq k \leq n$, $1 \leq i_1 < \dotsb < i_k \leq n$, $u(i_1) \dotsm u(i_k) s_{-1}$ maps $\ell^2(c+(d))$ bijectively into itself if and only if 
$c=-\sum_{m=1}^k d' \omega_{i_m}$, and a calculation shows that for this $c \in \cR_d$, $u(i_1) \dotsm u(i_k) s_{-1}$ acts on $\ell^2(c+(d))$ like $u(i_1)^* \dotsm u(i_k)^* s_{-1}$ under the identification $\ell^2(\so) \cong \ell^2(c+(d)); \xi_r \ma \xi_{c+dr}$. 

Thus, $\kappa_d(\eckl{\halb(1+u(i_1) \dotsm u(i_k) s_{-1})}_0) = \eckl{\halb(1+u(i_1) \dotsm u(i_k) s_{-1})}_0 + \tfrac{d^n-1}{2} \eckl{1}_0$. Then $\kappa_d$ is given by
\bgloz
  \rukl{
  \begin{array}{ccc|ccc}
  1 & & & & & \\
   & \ddots & & & 0 & \\
   & & 1 & & & \\
  \hline
  \tfrac{d^n-1}{2} & \dotso & \tfrac{d^n-1}{2} & d^n & & \\
   & 0 & & & \ddots & \\
   & & & & & d^{?}
\end{array}
}
\egloz
where $?$ is $0$ or $1$ depending on the parity of $n$.

Therefore, $\kappa_{2d}$ (for $d$ odd) is represented by the matrix
\bgloz
  \rukl{
  \begin{array}{cccc|ccc}
  1 & & & & & & \\
  \vdots & & 0 & & & 0 & \\
  1 & & & & & & \\
  \hline
  2^{n-1} d^n - 2^{n-1} & (2d)^{n-1} & \dotso & (2d)^{n-1} & (2d)^n & & \\
  \vdots & & 0 & & & \ddots & \\
  0 & & & & & & (2d)^{?}
\end{array}
}
\egloz
where $?$ is $0$ or $1$ depending on the parity of $n$. 

The result on the K-theory of $\fA^{(0)}$ follows.
\eproof

Our computations show that for any $d \in \Zz_{>1}$, $\Ad(s_d)$ induces $\id$ on the copies of $\Zz$ and multiplies the generators of $\Qz$ by some constant $>1$ in $K_0$. Thus, 
\bgl
\label{KA^0+d}
  K_j(C^*(\gekl{u^a}, s_{-1}, \gekl{e_b}, s_d)) 
  \cong 
  \bfa
    \Zz \falls n \text{ odd} \\
    \Zz^2 \fuer n \text{ even}
  \efa
\egl
for $j = 0, 1$. Here we have used the Pimsner-Voiculescu sequence together with the description $C^*(\gekl{u^a}, s_{-1}, \gekl{e_b}, s_d) \sim_MN \ilim \gekl{\fA^{(0)}; \Ad(s_d)} \rtimes_{\Ad(s_d)} \Zz$.

\subsection{Comparison between infinite and finite places}

The main result (Theorem \ref{M}) has already been proven in Section \ref{Me}. Applying Theorem \ref{M} to our situation, we obtain, together with Lemma \ref{KA^0}, 
\bgloz
  K_0(C_0(\Az_{\infty}) \rtimes K \rtimes \mu) 
  \cong 
  \bfa
    \Zz \oplus \Qz^{2^{n-1}} \falls n \text{ odd} \\
    \Zz \oplus \Qz^{2^{n-1}-1} \oplus \Zz \falls n \text{ even}, 
  \efa
\egloz
\bgloz
  K_1(C_0(\Az_{\infty}) \rtimes K \rtimes \mu) \cong \gekl{0} 
\egloz
and, using (\ref{KA^0+d}), 
\bgloz
  K_j(C_0(\Az_{\infty}) \rtimes K \rtimes (\mu \times \spkl{d})) 
  \cong 
  \bfa
    \Zz \falls n \text{ odd} \\
    \Zz^2 \falls n \text{ even} 
  \efa
\egloz
for $j=0,1$. Furthermore, studying the Pimsner-Voiculescu sequence associated to 
\bgloz
  C_0(\Az_{\infty}) \rtimes K \rtimes (\mu \times \spkl{d}) \cong (C_0(\Az_{\infty}) \rtimes K \rtimes \mu) \times \spkl{d}
\egloz
($\spkl{d} \cong \Zz$), we get that the only $\Qz \pos$-invariant part of $K_0(C_0(\Az_{\infty}) \rtimes K \rtimes \mu)$ is given by the copies of $\Zz$. We will need this observation later on.

\subsection{The multiplicative action}
\label{onlymu}

By equivariant Bott periodicity, we obtain 
\bgloz
  K_0(C_0(\Az_{\infty}) \rtimes \mu)
  \cong 
  \bfa
    \Zz \falls n \text{ odd} \\
    \Zz^2 \falls n \text{ even}  
  \efa
  \text{ and }
  K_1(C_0(\Az_{\infty}) \rtimes \mu) \cong \gekl{0}
\egloz
since $\Az_{\infty} \cong \Rz^n$ (compare Section \ref{ANT}).

Now, consider $i: C_0(\Az_{\infty}) \rtimes \mu \lori C_0(\Az_{\infty}) \rtimes K \rtimes \mu$, the homomorphism induced by the inclusion 
$C_0(\Az_{\infty}) \into C_0(\Az_{\infty}) \rtimes K$. Our aim is to prove that $i_*$ is injective on K-theory.

First of all, note that it suffices to prove that $i': C_0(\Az_{\infty}) \rtimes \mu \lori C_0(\Az_{\infty}) \rtimes \so \rtimes \mu$ is injective on $K_0$ since $C_0(\Az_{\infty}) \rtimes K \rtimes \mu$ can be written as an inductive limit (with $C_0(\Az_{\infty}) \rtimes \so \rtimes \mu$ as C*-algebra in each step) where the structure maps leave the image of $C_0(\Az_{\infty}) \rtimes \mu$ fixed on $K_0$. Thus, we have to prove the following

\blemma
\label{incinj}
$(i')_*$ is injective on $K_0$.
\elemma

\bproof
The main ingredient is Theorem A1 of \cite{Nat}, which is some sort of Mayer-Vietoris sequence relating the K-theory of crossed products by a free product of two groups with the K-theory of the single crossed products. We apply this result to the group $\so \rtimes \mu \cong \Zz^n \rtimes \Zz / 2 \Zz$. Namely, the identification $\Zz \rtimes  \Zz / 2 \Zz \cong \Zz / 2 \Zz * \Zz / 2 \Zz; z \ma t_2 t_1, t \ma t_1$ ($z, t, t_1, t_2$ are the canonical generators) is compatible with the action on $\Zz^{n-1}$ and on $C_0(\Az_{\infty})$. Thus, Theorem A1 of \cite{Nat} yields the following sequence which is exact in the middle:
\bgln
\label{nates}
  && K_0(C_0(\Az_{\infty}) \rtimes \Zz^k) \\
  &\overset{\kappa_*-\kappa'_*}{\lori}& 
  K_0(C_0(\Az_{\infty}) \rtimes \Zz^k \rtimes_{\hat{\beta}_{-1}} \Zz / 2 \Zz) 
  \oplus
  K_0(C_0(\Az_{\infty}) \rtimes \Zz^k \rtimes_{\hat{\alpha}_{\omega_{k+1}} \hat{\beta}_{-1}} \Zz / 2 \Zz) \nonumber \\
  &\overset{\ve_*+\ve'_*}{\lori}& 
  K_0(C_0(\Az_{\infty}) \rtimes \Zz^{k+1} \rtimes_{\hat{\beta}_{-1}} \Zz / 2 \Zz) \nonumber
\egln
where $\kappa$, $\kappa'$, $\ve$ and $\ve'$ are the canonical maps and $0 \leq k \leq n-1$.

Now, consider the translation-invariant isomorphism
\bgloz
  C_0(\Az_{\infty}) \lori C_0(\Az_{\infty}); f \ma f(\sqcup -\halb \omega_{k+1}).
\egloz
It yields an isomorphism
\bgloz
  C_0(\Az_{\infty}) \rtimes \Zz^k \overset{\psi}{\lori} C_0(\Az_{\infty}) \rtimes \Zz^k; f \cdot u^a \ma f(\sqcup -\halb \omega_{k+1}) \cdot u^a. 
\egloz
We have
\bglnoz
  && \psi \circ \hat{\beta}_{-1} (f \cdot u^a) = f(- \sqcup + \halb \omega_{k+1}) \cdot u^{-a} \\
  &=& \hat{\alpha}_{\omega_{k+1}} \circ \hat{\beta}_{-1} (f(\sqcup - \halb \omega_{k+1}) \cdot u^a) 
  = \hat{\alpha}_{\omega_{k+1}} \circ \hat{\beta}_{-1} (\psi(f \cdot u^a)).
\eglnoz
Thus, $\psi$ induces an isomorphism
\bgloz
  \psi: C_0(\Az_{\infty}) \rtimes \Zz^k \rtimes_{\hat{\beta}_{-1}} \Zz / 2 \Zz 
  \cong C_0(\Az_{\infty}) \rtimes \Zz^k \rtimes_{\hat{\alpha}_{\omega_{k+1}} \hat{\beta}_{-1}} \Zz / 2 \Zz.
\egloz
Now, the crucial point is that the following diagram commutes:
\bgloz
\begin{CD}
  C_0(\Az_{\infty}) \rtimes \Zz^k @> \kappa >> C_0(\Az_{\infty}) \rtimes \Zz^k \rtimes_{\hat{\beta}_{-1}} \Zz / 2 \Zz \\
  @VV \psi V @VV \psi V \\
  C_0(\Az_{\infty}) \rtimes \Zz^k @> \kappa' >> C_0(\Az_{\infty}) \rtimes \Zz^k \rtimes_{\hat{\alpha}_{\omega_{k+1}} \hat{\beta}_{-1}} \Zz / 2 \Zz  
\end{CD}
\egloz
This fact, together with $\psi \sim_h \id$ on $C_0(\Az_{\infty}) \rtimes \Zz^k$, implies that 
\bgl
\label{psi}
  \psi_* \kappa_* = \kappa'_*
\egl
on K-theory.

We would like to show that $\ve_*$ is injective. Assume that $x \in \ker(\ve_*)$. Then, $(x,0) \in \ker(\ve_* + \ve'_*) = \img(\kappa_* - \kappa'_*)$ because of (\ref{nates}). Thus, there exists $y$ in $K_0(C_0(\Az_{\infty}) \rtimes \Zz^k)$ with $\kappa_*(y) = x$ and $\kappa'_*(y) = 0$. But by (\ref{psi}), we have $0 = \kappa'_*(y) = \psi_* \kappa_*(y) = \psi_*(x)$, and since $\psi$ is an isomorphism, this implies $x = 0$. Thus, for any $0 \leq k \leq n-1$, $\ve_*$ is injective. This proves our claim, since $i' = \ve^{(n-1)} \circ \dotsb \circ \ve^{(0)}$ if $\ve^{(k)}$ denotes the $k$-th map 
\bgloz
  C_0(\Az_{\infty}) \rtimes \Zz^k \rtimes \mu \lori C_0(\Az_{\infty}) \rtimes \Zz^{k+1} \rtimes \mu.
\egloz
\eproof

\subsection{The general result}

It remains to put everything together. We distinguish between three cases:

a) $\# \gekl{v_{\Rz}} = 0$: 

Choose an unramified prime $p \in \Zz \pos$ and a $\Zz$-basis $\gekl{p, p_1, p_2, \dotsc}$ of $\Gamma$, where $K \reg = \mu \times \Gamma$. 
Let $\Gamma_m = \spkl{p, \dotsc, p_m}$. 

We have seen $K_j(C_0(\Az_{\infty}) \rtimes K \rtimes (\mu \times \Gamma_0)) \cong \Zz^2$ for $j=0,1$ ($\# \gekl{v_{\Rz}} = 0$ implies that $n$ is even). Moreover, it follows from our results of Section \ref{onlymu} and the Pimsner-Voiculescu sequence that 
$i_*: K_*(C_0(\Az_{\infty}) \rtimes (\mu \times \Gamma_0)) \lori K_*(C_0(\Az_{\infty}) \rtimes K \rtimes (\mu \times \Gamma_0))$ 
is injective. Now, let us prove inductively that 
\bgloz
  K_j(C_0(\Az_{\infty}) \rtimes (\mu \times \Gamma_m)) \cong K_j(C_0(\Az_{\infty}) \rtimes K \rtimes (\mu \times \Gamma_m)) \cong \Zz^{2^{m+1}}
\egloz
for $j=0,1$ and that $i_*: K_*(C_0(\Az_{\infty}) \rtimes (\mu \times \Gamma_m)) \lori K_*(C_0(\Az_{\infty}) \rtimes K \rtimes (\mu \times \Gamma_m))$ is injective. The case $m=0$ has already been proven above. 

If the claim is proven for $m$, the Pimsner-Voiculescu sequence, together with its naturality, will yield the result for $m+1$ (analogously to the first case of Section \ref{KB}, we have to use that $\hat{\beta}_b \sim_h \id$ equivariantly).

Thus, we conclude $K_*(\fA) \cong K_*(C_0(\Az_{\infty}) \rtimes K \rtimes K \reg) \cong K_0(C^*(\mu)) \otimes_{\Zz} \extalg (\Gamma)$.

b) $\# \gekl{v_{\Rz}}$ odd: 

Again, choose $p \in \Zz \pos$ prime and unramified and a $\Zz$-basis $\gekl{p, p_1, \dotsc}$ of $\Gamma$ with $K \reg = \mu \times \Gamma$. As $\# \gekl{v_{\Rz}}$ is odd, we can arrange by multiplying with $-1$ that $\# \menge{v_{\Rz}}{v_{\Rz}(p_i)<0}$ is even for all $i$. As above, let $\Gamma_m = \spkl{p, \dotsc, p_m}$. We can show for each $m$ that 
\bgloz
  K_j(C_0(\Az_{\infty}) \rtimes (\mu \times \Gamma_m)) \cong K_j(C_0(\Az_{\infty}) \rtimes K \rtimes (\mu \times \Gamma_m)) \cong \Zz^{2^m}
\egloz
for $j=0,1$ and that $i_*: K_*(C_0(\Az_{\infty}) \rtimes (\mu \times \Gamma_m)) \lori K_*(C_0(\Az_{\infty}) \rtimes K \rtimes (\mu \times \Gamma_m))$ is injective. 

All we have to show is that $(\hat{\beta}_{p_{m+1}})_* = \id$ on $K_*(C_0(\Az_{\infty}) \rtimes (\mu \times \Gamma_m))$. This follows from $\hat{\beta}_{p_{m+1}} \sim_h \id$ on $C_0(\Az_{\infty}) \rtimes \mu$ and that $\hat{\beta}_{p_{m+1}} \sim_h \hat{\beta}_{(v(p_{m+1}) / \abs{v(p_{m+1})})_{v \tei \infty}}$. The second fact implies that $\hat{\beta}_{p_{m+1}}$ is of period $2$, while the first one, together with the Pimsner-Voiculescu sequence, shows that $\hat{\beta}_{p_{m+1}}$ can be described by an upper triangular matrix where all the diagonal elements are $1$. These two facts imply our claim, namely $(\hat{\beta}_{p_{m+1}})_* = \id$ on $K_*(C_0(\Az_{\infty}) \rtimes (\mu \times \Gamma_m))$. 

Thus, $K_*(\fA) \cong K_*(C_0(\Az_{\infty}) \rtimes K \rtimes K \reg) \cong \extalg (\Gamma)$.

c) $\# \gekl{v_{\Rz}} \geq 2$ even: 

Again, let $K \reg = \mu \times \Gamma$ and choose a $\Zz$-basis $\gekl{p, p_1, p_2, \dotsc}$ of $\Gamma$, with $p \in \Zz \pos$. We can arrange that 
$\# \menge{v_{\Rz}}{v_{\Rz}(p_1)<0}$ is odd and $\# \menge{v_{\Rz}}{v_{\Rz}(p_i)<0}$ is even for all $i>1$. Let $\Gamma_m = \spkl{p, \dotsc, p_m}$ and $\Gamma'_m = \spkl{p, p_2 \dotsc, p_m}$. As above, we can show that $K_j(C_0(\Az_{\infty}) \rtimes K \rtimes (\mu \times \Gamma'_m)) \cong \Zz^{2^m}$ for $j=0,1$ ($\eckl{K:\Qz}$ even by assumption). Additionally, we have 

\blemma
\label{1-1}
$\hat{\beta}_{p_1} 
  = 
  \rukl{
  \begin{smallmatrix}
  1 & & & & & & \\
   & -1 & & & & 0 & \\
   & & 1 & & & & \\
   & & & -1 & & & \\
   & & &  & \ddots & & \\
   & 0 & & & & 1 & \\
   & & & & & & -1 \\
  \end{smallmatrix}
  }
$
on $K_j(C_0(\Az_{\infty}) \rtimes K \rtimes (\mu \times \Gamma'_m)) \cong \Zz^{2^m}$ for $j=0,1$ with respect to an appropriate $\Zz$-basis.
\elemma

\bproof
$\hat{\beta}_{p_1} \sim_h \hat{\beta}_{(v(p_1) / \abs{v(p_1)})_{v \tei \infty}}$ and $\hat{\beta}_{(v(p_1) / \abs{v(p_1)})_{v \tei \infty}} = \hat{\beta}_{(1,-1,1,1, \dotsc) \otimes 1_{\Cz}} \circ \hat{\beta}_{p'}$ where $\# \menge{v_{\Rz}}{v_{\Rz}(p')<0}$ is even. Thus, $(\hat{\beta}_{p'})_* = \id$ on $K_*(C_0(\Az_{\infty}) \rtimes K \rtimes (\mu \times \Gamma'_m))$ and it remains to show the claim for $\hat{\beta}_{(1,-1,1,1, \dotsc) \otimes 1_{\Cz}})$. We proceed inductively: 

To get started, consider the special case $K = \Qz \eckl{\sqrt{2}}$. Then $\so = \Zz + \Zz \sqrt{2}$. We would like to show 
$(\hat{\beta}_{\sqrt{2}})_* = 
  \rukl{
  \begin{smallmatrix}
  1 & \\
   & -1 
  \end{smallmatrix}
  }
$ 
on $K_0(C_0(\Az_{\infty}) \rtimes \mu)$. By our Theorem \ref{M} and the Pimsner-Voiculescu sequence applied to $C_0(\Az_{\infty}) \rtimes (\mu \times \spkl{p})$, it suffices to show that $\Ad(s_{\sqrt{2}})$ induces 
$
  \rukl{
  \begin{smallmatrix}
  1 & \\
   & -1 
  \end{smallmatrix}
  }
$
on $K_j(C^*(\gekl{u^a}, s_{-1}, \gekl{e_b}, s_p))$ ($j=0,1$).

To determine $\Ad(s_{\sqrt{2}})$, we compute $(e_{0,0} \otimes \id)_*^{-1} (\varphi_{\sqrt{2}})_* (\iota_{\sqrt{2}})_*$. We have
\bgloz
  (\nu_{\sqrt{2}})_*(\eckl{u^1}_1 \times \eckl{u^{\sqrt{2}}}_1) = \eckl{u^{\sqrt{2}}}_1 \times \eckl{u^2}_1 = -2 \eckl{u^1}_1 \times \eckl{u^{\sqrt{2}}}_1.
\egloz
Thus $\Ad(s_{\sqrt{2}})$ induces $-1$ on the second copy of $\Zz$ in $K_0(C^*(\gekl{u^a}, s_{-1}, \gekl{e_b}))$. Here we used that (\ref{cru}) also holds for $d = \sqrt{2}$. 

Now, take $\cR_{\sqrt{2}} = \gekl{0, 1}$. Then, under similar identifications as in Lemma \ref{mat}, 
\bgloz
  (\varphi_{\sqrt{2}}) \circ (\iota_{\sqrt{2}}) (\halb(1+s_{-1})) = 
  \rukl{
  \begin{smallmatrix}
  \halb(1+s_{-1}) & \\
   & \halb(1+u^{-\sqrt{2}}s_{-1})
  \end{smallmatrix}
  }
\egloz 
and 
\bgloz
  (\varphi_{\sqrt{2}}) \circ (\iota_{\sqrt{2}}) (\halb(1+u^1 s_{-1})) = 
  \halb
  \rukl{
  \begin{smallmatrix}
  1 & s_{-1} \\
  s_{-1} & 1
  \end{smallmatrix}
  }.
\egloz 
Thus, 
\bglnoz
  && (e_{0,0} \otimes \id)_*^{-1} \circ (\varphi_{\sqrt{2}})_* \circ (\iota_{\sqrt{2}})_* (\eckl{\halb(1+s_{-1})}_0 - \eckl{\halb(1+u^1 s_{-1})}_0) \\
  &=& \eckl{\halb(1+s_{-1})}_0 -(\eckl{1}_0 - \eckl{\halb(1+u^{-\sqrt{2}} s_{-1})}_0).
\eglnoz
In the inductive limit, we get that $\Ad(s_{\sqrt{2}})$ induces $\id_{\Zz}$ on the first copy of $\Zz$ in $K_0(C^*(\gekl{u^a}, s_{-1}, \gekl{e_b}))$ by the analogue of (\ref{KRel}). Thus, $(\hat{\beta}_{(1,-1)})_* = 
  \rukl{
  \begin{smallmatrix}
  1 & \\
   & -1 
  \end{smallmatrix}
  }
$ on $K_0(C_0(\Az_{\infty}) \rtimes \mu) = K_0(C_0(\Rz^2) \rtimes \mu)$ as $K \ni \sqrt{2} \ma (\sqrt{2},-\sqrt{2}) \sim_h (1,-1) \in \Az_{\infty}$. 

But this already implies $(\hat{\beta}_{(1,-1,1,1, \dotsc) \otimes 1_{\Cz}})_* 
  = 
  \rukl{
  \begin{smallmatrix}
  1 & \\
   & -1 
  \end{smallmatrix}
  }
$ on $K_0(C_0(\Az_{\infty}) \rtimes \mu)$ for any number field with $\mu = \gekl{\pm 1}$: Let $B$ be the Bott element in $KK^{\mu}(\Cz, C_0(\Rz^{n-2}))$. Then $1_{C_0(\Rz^2)} \otimes B$ is an invertible element of $KK^{\mu}(C_0(\Rz^2), C_0(\Rz^2) \otimes C_0(\Rz^{n-2}))$. By \cite{Kas}, Theorem 2.14. 8),  
\bgloz
  (\hat{\beta}_{(1,-1)})_* \times (1_{C_0(\Rz^2)} \otimes B) \times \sigma_* 
  = (B \otimes 1_{C_0(\Rz^2)}) \times (1_{C_0(\Rz^{n-2})} \otimes (\hat{\beta}_{(1,-1)})_*) 
\egloz
where $\sigma$ is the flip $C_0(\Rz^2) \otimes C_0(\Rz^{n-2}) \cong C_0(\Rz^{n-2}) \otimes C_0(\Rz^2)$. This implies
\bglnoz
  && (\hat{\beta}_{(1,-1)})_* \times (1_{C_0(\Rz^2)} \otimes B) \\
  &=& (B \otimes 1_{C_0(\Rz^2)}) \times (1_{C_0(\Rz^{n-2})} \otimes (\hat{\beta}_{(1,-1)})_*) \times \sigma_*^{-1} \\
  &=& (B \otimes 1_{C_0(\Rz^2)}) \times \sigma_*^{-1} \times ((\hat{\beta}_{(1,-1)})_* \otimes 1_{C_0(\Rz^{n-2})}) \\
  &=& (1_{C_0(\Rz^2)} \otimes B) \times ((\hat{\beta}_{(1,-1)})_* \otimes 1_{C_0(\Rz^{n-2})}).
\eglnoz
Thus, using $1_{C_0(\Rz^2)} \otimes B$ to identify $K_0(C_0(\Az_{\infty}) \rtimes \mu)$ with $\Zz^2 \cong K_0(C_0(\Rz^2) \rtimes \mu)$, 
$(\hat{\beta}_{(1,-1,1,1, \dotsc) \otimes 1_{\Cz}})_*$ is given by 
$ 
  \rukl{
  \begin{smallmatrix}
  1 & \\
   & -1 
  \end{smallmatrix}
  }
$. This proves our claim for $m=1$. 

To go from $m$ to $m+1$, we apply the induction hypothesis together with the Pimsner-Voiculescu sequence to find a $\Zz$-basis for $K_j(C_0(\Az_{\infty}) \rtimes K \rtimes (\mu \times \Gamma'_{m+1}))$ for $j=0,1$ such that 
\bgloz 
  (\hat{\beta}_{(1,-1,1,1, \dotsc) \otimes 1_{\Cz}})_* 
  =
  \rukl{
  \begin{array}{ccccc|ccccc}
   1 & & & & & & & & & \\
   & -1 & & 0 & & & & & & \\
   & & \ddots & & & & & * & & \\
   & 0 & & 1 & & & & & & \\
   & & & & -1 & & & & & \\
  \hline
   & & & & & 1 & & & & \\
   & & & & & & -1 & & 0 & \\
   & & 0 & & & & & \ddots & & \\
   & & & & & & 0 & & 1 & \\
   & & & & & & & & & -1 \\
  \end{array}
  }.
\egloz
Now, the Pimsner-Voiculescu sequence implies that the torsion-free part of $K_j(C_0(\Az_{\infty}) \rtimes K \rtimes (\mu \times \Gamma_{m+1}))$ is $\Zz^{2^{m+1}}$. Thus, using $(\hat{\beta}_{(1,-1,1,1, \dotsc) \otimes 1_{\Cz}})_*^2 = \id$, we deduce that 
$K_j(C_0(\Az_{\infty}) \rtimes K \rtimes (\mu \times \Gamma'_{m+1})) / \img(\id-(\hat{\beta}_{(1,-1,1,1, \dotsc) \otimes 1_{\Cz}})_*)$ can only contain torsion of order $2$. Therefore the Pimsner-Voiculescu sequence implies $K_j(C_0(\Az_{\infty}) \rtimes K \rtimes (\mu \times \Gamma'_{m+1})) / \img(\id-(\hat{\beta}_{(1,-1,1,1, \dotsc) \otimes 1_{\Cz}})_*) \cong \Zz^{2^m} \oplus (\Zz / 2 \Zz)^{2^m}$. But this result, together with the elementary divisor theorem, tells us that we can modify the first chosen $\Zz$-basis for $K_j(C_0(\Az_{\infty}) \rtimes K \rtimes (\mu \times \Gamma'_{m+1}))$ so that our claim holds.
\eproof

Hence, applying the Pimsner-Voiculescu sequence iteratively gives 
\bgloz
  K_*(C_0(\Az_{\infty}) \rtimes K \rtimes (\mu \times \Gamma_m)) \cong \extalg (\Gamma'_m) \oplus (\Zz / 2 \Zz \otimes_{\Zz} \extalg (\Gamma'_m))
\egloz
and thus 
\bgloz
  K_*(\fA) \cong \extalg (\Gamma) \oplus ((\Zz / 2 \Zz) \otimes_{\Zz} \extalg (\Gamma)).
\egloz

\bremark
At this point, it also becomes clear that we cannot treat the general case (without the restriction $\mu = \gekl{\pm 1}$) for the following reasons: It is not clear how to prove analogous statements as Lemma \ref{BC} and Lemma \ref{incinj} in general. Once these two problems are solved, it should be possible to determine the K-theory without further assumptions on the number fields. 
\eremark

\bremark
With a similar idea as in Theorem \ref{M}, we can treat the case of the full adele ring with the action of the $ax+b$-group: Let $K$ be any global field. As $K$ is a discrete subgroup of $\Az$, it acts freely and properly on $\Az$ so that 
\bgloz
  C_0(\Az) \rtimes K \sim_M C(\Az / K) \cong C(\hat{K}) \cong C^*(K).
\egloz
$C^*(K)$ is the group C*-algebra of $(K,+)$. Moreover, it turns out that this Morita equivalence can be chosen equivariantly (in the sense of \cite{CMW}) with respect to the multiplicative action of $K \reg$. Thus, by \cite{CMW}, we get
\bgl
\label{Mfull}
  C_0(\Az) \rtimes P_K = C_0(\Az) \rtimes K \rtimes K \reg \sim_M C^*(K \rtimes K \reg) = C^*(P_K).
\egl
This means that the crossed product is Morita equivalent to the C*-algebra of the $ax+b$-group over $K$. We note that this group C*-algebra is also the ring C*-algebra associated with the field $K$.

(\ref{Mfull}) can be used to compute the K-theory. For example, in the case $K = \Qz$, we get 
\bgloz
  K_*(C_0(\Az) \rtimes \Qz \rtimes \Qz \reg) \cong K_0(C^*(\gekl{\pm 1})) \otimes_{\Zz} \extalg (\Qz \pos).
\egloz
\eremark

\end{document}